\documentclass[11pt]{article}
\usepackage[tbtags]{amsmath}
\usepackage{amssymb}
\usepackage{amsthm}
\usepackage[misc]{ifsym}
\usepackage{cases}
\usepackage{cite}	
\usepackage{mathrsfs}	
\usepackage{xcolor}	
\usepackage{graphicx}	
\usepackage{float}	
\usepackage{color}	
\usepackage{hyperref}	

\hypersetup{hypertex=true,
	colorlinks=true,
	linkcolor=red,
	anchorcolor=blue,
	citecolor=blue}

\numberwithin{equation}{section}
\normalsize

\title{\bf Linear-Quadratic Non-zero Sum Differential Game with Asymmetric Delayed Information
	\thanks{This work is financially supported by the National Key R\&D Program of China (2022YFA1006104), National Natural Science Foundations of China (12471419, 12271304), and Shandong Provincial Natural Science Foundations (ZR2024ZD35, ZR2022JQ01).}}
\author{\normalsize
	Yuxin Ye\thanks{\textit{School of Mathematics, Shandong University, Jinan 250100, P.R.  China, E-mail: yeyuxin@mail.sdu.edu.cn}},
	\ Jingtao Shi\thanks{\textit{Corresponding author, School of Mathematics, Shandong University, Jinan 250100, P.R. China, E-mail: shijingtao@sdu.edu.cn}}}	
\date{}

\newtheorem{mythm}{Theorem}[section]

\newtheorem{mylem}{Lemma}[section]
\newtheorem{Remark}{Remark}[section]

\begin{document}
	
\maketitle

\noindent{\bf Abstract:}\quad
	 This paper is concerned with a  linear-quadratic non-zero sum differential game with asymmetric delayed information. To be specific, two players exist time delays simultaneously which are different, leading the dynamical system being an asymmetric information structure. By virtue of stochastic maximum principle,  the stochastic Hamiltonian system is given which is a delayed forward-backward stochastic differential equation. Utilizing discretisation approach and backward iteration technique, we establish the relationship between forward and backward processes under asymmetric delayed information structure and obtain the state-estimate feedback Nash equilibrium of our problem.
	 
	\vspace{2mm} 
	
	
	\vspace{2mm}
	
	
\section{Introduction}

  	Non-zero sum stochastic differential game is a significant type of differential games, and the pioneering contributions can be traced back to Starr and Ho \cite{Starr1969}, in which they raised the mathematical formulation of this type of game and analyzed its solution. In non-zero sum differential games, which are not completely competitive, all players have the same objective, that is, they would like to maximize their own profits or minimize their costs through joint control. This optimal control strategy is called Nash equilibrium. It has extensive applications in finance and insurance (Bensoussan et al. \cite{BS2014}). There are also many works on the {\it linear-quadratic} (LQ) case. Hamadène \cite{Hamadene1999} and Wu \cite{Wu2005} studied the LQ non-zero sum stochastic differential games and built the relationship between the solvability of the corresponding {\it forward-backward stochastic differential equations} (FBSDEs) and the existence of Nash equilibria. Sun and Yong \cite{Sun2019} explored the LQ non-zero sum stochastic differential game under the concepts of open-loop and closed-loop Nash equilibria. More research works on non-zero sum stochastic differential game can be found in Wu and Yu \cite{WY2005}, Wang and Yu \cite{WY2010, WY2012}, Wang et al. \cite{WXX2018}, Nie et al.\cite{NWY2022} and the references therein.
  	
  	System with time delay is such a system that the performance of the system at a moment depends on the history in the past. In recent years, delayed system has attracted the attention of researchers due to its wide applications in such as optimal consumption problems in the economy (Agram and \O ksendal \cite{AO2014}), finance (Arriojas et al. \cite{Arriojas2007}), wireless communication network (Bansal and Liu \cite{BL2003}) an so on. Chen and Wu \cite{CW2010} obtained the maximum principle for the delayed stochastic control problem utilizing classical variational method and duality method. Agram and \O ksendlal \cite{AO2014} studied the maximum principle for the infinite-dimensional optimal control problem under FBSDEs with time delay. Meng et al. \cite{Meng2025} extended the stochastic maximum principle of the delayed system using the Volterra integral equation. Shi and Wang \cite{SW2015} studied the theory and application of non-zero sum differential game with state equation being a delayed BSDEs. There are many other research works on stochastic control problems of delayed system, such as Chen et al. \cite{CWY2012}, Zhang and Xu \cite{ZX2017}, Wang et al.\cite{WZX2021} and the references therein.
  	
  	It is worth noting that special FBSDEs are involved in the stochastic Hamiltonian system of delayed control problems. In Peng and Yang \cite{PY2009}, the solvability of {\it anticipated BSDEs} (ABSDEs) was studied, in which the conditional expectation of future state was contained. Meanwhile, Xu and Zhang \cite{XZ2018}, Xu et al. \cite{XZX2018}, Ma et al. \cite{Ma2022, Ma2024} discussed the explicit solutions of {\it delayed FBSDEs} (DFBSDEs) by introducing modified Riccati equations, and using discretisation approach and backward iteration technique.
  	
  	In this paper, we explore the LQ non-zero sum stochastic differential game with different delays, which form an asymmetric information structure. At present, there are many studies based on asymmetric information. Shi et al. \cite{SWX2016, SWX2017} investigated the Stackelberg stochastic differential game under asymmetric information, where the information filtrations of two player are sub-$\sigma$-algebra of the general filtration $\mathcal{F}_t$. Wang et al. \cite{WXX2018} discussed the non-zero sum game under asymmetric information with the state equation being a BSDE. Ma et al. \cite{MXW2022} studied an LQ optimal control problem for the stochastic system under asymmetric information structure. Chen and Zhang \cite{CZ2026} investigated one kind of partial information non-zero sum stochastic differential game with mixed delays. In this paper, we employ different method to search for the Nash equilibrium, and we also establish the relationship between the forward and backward processes. Firstly, by leveraging the maximum principle of stochastic system with time delay \cite{CW2010}, we transform the solvability of the original problem into that of DFBSDEs with two different conditional expectations. Secondly, we adopt the discretisation approach and the backward iteration technique in \cite{Ma2024} to decouple the DFBSDEs, thereby we obtain the explicit solutions in discrete-time form. By taking the limit of the time interval, we formally derive the explicit solutions in continuous-time form and the corresponding Riccati equations, and rigorously verify that they are indeed the solutions of initial DFBSDEs. Finally, we present the state-estimate feedback Nash equilibrium of the original problem.
  	
  	The rest of the paper is organized as follows. Problem formulation is stated in Section 2. In Section 3, we solve the problem and derive the state-estimate feedback Nash equilibrium. Conclusion is provided in Section 4.
  	
 	\textit{Notations}: $\mathbb{R}^n$ stands for the $n$-dimensional Euclidean space. The superscript $'$ stands for transpose of a matrix. \textit{I} represents the identity matrix of appropriate dimensions. Let $[0, T ]$ be a finite time horizon and $(\Omega,\mathcal{F},\{\mathcal{F}_t\}_{t\geq0},\mathbf{P})$ be a complete filtered probability space with the filtration defined as $\{\mathcal{F}_t\}_{t\geq0}:=\sigma\{w(s);0 \leq s\leq t\}\}\bigvee\mathcal{N}$, where $w(\cdot)$ is a $\mathit{d}$-dimensional standard Brownian motion and $\mathcal{N}$ contains all $\mathbf{P}$-null sets in $\mathcal{F}$. $\mathbb{E}_s [x(t)]:=\mathbb{E}[x(t)|\mathcal{F}_s] $ denotes the conditional expectation of $x(t)$ with regard to $\mathcal{F}_s$, $0\leq s\leq t\leq T$.

\section{Problem formulation}

	In this paper, we consider the following controlled linear SDE on $[0, T]$:
	\begin{equation}
		\begin{cases}\label{linear SDE}
	dx(t)=[Ax(t)+B_1v_1(t)+B_2v_2(t)]dt+[\bar{A}x(t)+\bar{B}_1v_1(t)+\bar{B}_2v_2(t)]dw(t),\\x(0)=x_0,
	\end{cases}
	\end{equation} 
	where $x(\cdot)\in \mathbb{R}^n$ is the system state, $v_i(\cdot)\in \mathbb{R}^{d_i}$ is the control process of player $i, i=1,2$. $A, B_1, B_2, \bar{A}, \bar{B}_1, \bar{B}_2 $ are constant matrices with compatible dimensions. The initial value $x_0$ is given in advance and independent of $w(\cdot)$.
	
	Now we would like to explain the asymmetric delayed information feature in this paper. For two players, there exist two different time delays $h_i (i=1,2,  0<h_2<h_1<T)$. The information available to player 1 at time $t$ is based on $\mathcal{F}_{t-h_1}$, while $\mathcal{F}_{t-h_2}$ is the information accessed by player 2, thus $\mathcal{F}_{t-h_1} \subseteq \mathcal{F}_{t-h_2} \subseteq \mathcal{F}_t$. We denote by $\mathcal{U}_{ad}^i$ the set of all $\mathbb{R}^{d_i}$-valued $\mathcal{F}_{t-h_i}$-adapted square-integrable control processes $v_i(\cdot)$ (i.e., $ \mathbb{E}\int_0^T|v_i(t)|^2dt<\infty, i=1, 2$). Each element in $\mathcal{U}_{ad}^i[0, T]$ ia called an (open-loop) admissible control of player $i$, $i=1,2$. 
	
	In non-zero sum stochastic differential game, each player has its own cost hoping to minimize, which can be described by the cost functionals as follows: 
	\begin{equation}\label{cost functionals}
		\begin{aligned}
			\mathit{J}_1(v_1(\cdot))=\frac{1}{2}\mathbb{E}\left\{\int_0^T[x'(t)Q_1x(t)+v'_1(t)R_1v'_1(t)]dt+x'(T)H_1x(T)\right\},\\
			\mathit{J}_2(v_2(\cdot))=\frac{1}{2}\mathbb{E}\left\{\int_0^T[x'(t)Q_2x(t)+v'_2(t)R_2v'_2(t)]dt+x'(T)H_2x(T)\right\}.
		\end{aligned}
	\end{equation} 
	Here $Q_i, H_i $ are $n\times n$-dimensional positive semi-definite matrices, $R_i$ is $d_i \times d_i$-dimensional positive definite matrix, $i=1,2$. The LQ non-zero sum stochastic differential game with asymmetric delayed information is stated as the following.
	
	\noindent\textbf{Problem SDG-ADI:} Find  $(u_1(\cdot), u_2(\cdot))\in \mathcal{U}_{ad}^1 \times\mathcal{U}_{ad}^2$ such that
	\begin{equation}
		\begin{cases}
			\mathit{J}_1(u_1(\cdot))\leq\mathit{J}_1(v_1(\cdot)), \quad\forall v_1(\cdot)\in \mathcal{U}_{ad}^1, \\
			\mathit{J}_2(u_2(\cdot))\leq \mathit{J}_2(v_2(\cdot)),\quad\forall v_2(\cdot)\in \mathcal{U}_{ad}^2,
		\end{cases}
	\end{equation} 
	which is subject to (\ref{linear SDE}).
	
	\begin{Remark}
	Due to the existence of two different delays, the solvability of this problem will become much more complex than the case with one delay.
	\end{Remark}

\section{Main results}
	
	The LQ non-zero sum stochastic differential game described in the previous section, can be viewed as two stochastic LQ optimal control problems with asymmetric delayed information structure. In order to obtain the Nash equilibrium, through stochastic maximum principle, the solvability of the problem is reduced to that of the corresponding DFBSDEs. Inspired by the works of Ma et al. \cite{Ma2024}, we obtain the explicit solution in discrete-time form by using backward iterative induction, and the continuous-time form solution by taking limit of the time interval. At last, the state-estimate feedback Nash equilibrium is presented.
		
	To begin with, by using stochastic maximum principle introduced in \cite{CW2010}, the solvability of {\bf Problem SDG-ADI} is given as follows.
		 
	\begin{mylem}
		 If $(u_1(\cdot),u_2(\cdot))\in \mathcal{U}_{ad}^1\times\mathcal{U}_{ad}^2$ is a Nash equilibrium for \textbf{Problem SDG-ADI}, then it satisfies the following equalities:
		  	\begin{align}
		  		0=R_1u_1(t)+\mathbb{E}_{t-h_1}[B'_1p_1(t)+\bar{B}'_1q_1(t)],\label{Nash equilibrium1} \\
                0=R_2u_2(t)+\mathbb{E}_{t-h_2}[B'_2p_2(t)+\bar{B}'_2q_2(t)],\label{Nash equilibrium2}
		  	\end{align}
		where $(p_i(\cdot),q_i(\cdot)), i=1,2$ satisfies:
		  	\begin{equation}
		  		\begin{cases}\label{adjoint equations}
		  			-dp_i(t)=[A'p_i(t)+\bar{A}'q_i(t)+Q_ix(t)]dt-q_i(t)dw(t),\quad t\in[0,T], \\p_i(t)=H_ix(T),
		  		\end{cases}
		  	\end{equation}
        and now $x(\cdot)\equiv x(\cdot;u_1(\cdot),u_2(\cdot))$ denoted the solution of (\ref{linear SDE}) corresponding to $(u_1(\cdot),u_2(\cdot))$.
		  	\end{mylem}
		
		Lemma 3.1 shows that if {\bf Problem SDG-ADI} is solvable, its Nash equilibrium $(u_1(\cdot),u_2(\cdot))$ can be characterized by the following FBSDEs:
			\begin{equation}
				\begin{cases}\label{FBSDEs}
					dx(t)=[Ax(t)+B_1u_1(t)+B_2u_2(t)]dt+[\bar{A}x(t)+\bar{B}_1u_1(t)+\bar{B}_2u_2(t)]dw(t), \\
                    -dp_i(t)=[A'p_i(t)+\bar{A}'q_i(t)+Q_ix(t)]dt-q_i(t)dw(t), \\
                    0=R_1u_1(t)+\mathbb{E}_{t-h_1}[B'_1p_1(t)+\bar{B}'_1q_1(t)], \\
                    0=R_2u_2(t)+\mathbb{E}_{t-h_2}[B'_2p_2(t)+\bar{B}'_2q_2(t)],\\
                    x(0)=x_0,\quad p_i(T)=H_ix(T),\quad i=1,2.
				\end{cases}
			\end{equation}
		
		Moreover, the solution of the FBSDEs (\ref{FBSDEs}) is unique if the optimal control is unique, which is the result below. 
		\begin{mylem}
			{\bf Problem SDG-ADI} is solvable if and only if there exits a unique $(x(\cdot),p_1(\cdot),q_1(\cdot),\\p_2(\cdot),q_2(\cdot))$ satisfying FBSDEs (\ref{FBSDEs}).
		\end{mylem}	  	 
		\noindent\textit{Proof.} As the proof is similar to those for Proposition 5.5, Chapter 6 in \cite{YZ1999} and Theorem 2.1 in \cite{CW2010}, it is omitted here.$\hfill\qedsymbol$

		Considering $R_1, R_2$ are positive definite matrices, thus (\ref{Nash equilibrium1}), (\ref{Nash equilibrium2}) can be transformed as follows:
			\begin{align}
				u_1(t)=-R_1^{-1}\mathbb{E}_{t-h_1}[B'_1p_1(t)+\bar{B}'_1q_1(t)], \\u_2(t)=-R_2^{-1}\mathbb{E}_{t-h_2}[B'_2p_2(t)+\bar{B}'_2q_2(t)].
			\end{align}
		
		Then FBSDEs (\ref{FBSDEs}) is equivalent to the following DFBSDEs:
			\begin{equation}
				\begin{cases}\label{DFBSDEs}
					dx(t)=\Big\{Ax(t)+B_{11}\mathbb{E}_{t-h_1}[p_1(t)]+B_{12}\mathbb{E}_{t-h_2}[p_2(t)]+B_{21}\mathbb{E}_{t-h_1}[q_1(t)]\\
					\qquad\qquad +B_{22}\mathbb{E}_{t-h_2}[q_2(t)]\Big\}dt
					+\Big\{\bar{A}x(t)+\bar{B}_{11}\mathbb{E}_{t-h_1}[p_1(t)]+\bar{B}_{12}\mathbb{E}_{t-h_2}[p_2(t)]\\
					\qquad\qquad +\bar{B}_{21}\mathbb{E}_{t-h_1}[q_1(t)]+\bar{B}_{22}\mathbb{E}_{t-h_2}[q_2(t)]\Big\}dw(t), \\
                    -dp_i(t)=[A'p_i(t)+\bar{A}'q_i(t)+Q_ix(t)]dt-q_i(t)dw(t), \\x(0)=x_0,\quad p_i(T)=H_iX(T),\quad i=1,2,
				\end{cases}
			\end{equation}
		where 
		\begin{align*}
			B_{11}:=-B_1R_1^{-1}B'_1,\quad B_{12}:=-B_2R_2^{-1}B'_2,\quad B_{21}:=-B_1R_1^{-1}\bar{B}'_1,\quad B_{22}:=-B_2R_2^{-1}\bar{B}'_2,\\ 
            \bar{B}_{11}:=-\bar{B}_1R_1^{-1}B'_1,\quad \bar{B}_{12}:=-\bar{B}_2R_2^{-1}B'_2,\quad \bar{B}_{21}:=-\bar{B}_1R_1^{-1}\bar{B}'_1,\quad \bar{B}_{22}:=-\bar{B}_2R_2^{-1}\bar{B}'_2.
		\end{align*}
		
		The solvability of {\bf Problem SDG-ADI} is now equivalent to that of the DFBSDEs (\ref{DFBSDEs}). In order to obtain the solution of (\ref{DFBSDEs}), we utilize the method of discretisation and the technique of backward iteration. Then by taking limit, we could access the continuous-time form solution.

		Before that, following the same procedure introduced in Ma et al. \cite{Ma2024}, we divide $[0,T ]$ into $N+1$ equal parts, i.e, $0=t_0<t_1<\dots<t_{N+1}=T$, and we choose an appropriate $\delta$ which is sufficiently small as the length of the subinterval, ensuring $h_i=d_i\delta$, and $d_i$ being integer, $i=1,2$, $0<d_2<d_1<N+1$. Denote $\Delta w_k:= w_{k+1}-w_k$ and $\hat{\mathcal{F}}_k$ is the natural filtration generated by $\Delta w_k, i.e., \hat{\mathcal{F}}_k=\sigma\{\Delta w_0,\dots, \Delta w_k\}$, where $\Delta w_k$ is a white noise with zero mean and covariance $\delta$. Simply denote the variables in (\ref{DFBSDEs}) at time $t_k$ as $(x_k, p_{i,k}, q_{i,k})$, $i=1,2$. We introduce the following notations: 
		\begin{align*}
			A'_k:= \hat{A}'+\Delta w_k\bar{A}',\quad A_k:= \hat{A}+\Delta w_k\bar{A},\quad \hat{Q}_{i}:=\delta Q,\\
            B_{11,k}:=\hat{B}_{11}+\Delta w_k\bar{B}_{11},\quad B_{12,k}:=\hat{B}_{12}+\Delta w_k\bar{B}_{12}, \\
            B_{21,k}:=\hat{B}_{21}+\Delta w_k\bar{B}_{21},\quad B_{22,k}=\hat{B}_{22}+\Delta w_k\bar{B}_{22}, 
		\end{align*}
		where $	\hat{A'}:= I+\delta A',\hat{A}:=I+\delta A, \hat{B}_{11}:=\delta B_{11},\hat{B}_{12}:=\delta B_{12}, \hat{B}_{21}:=\delta B_{21}, \hat{B}_{22}:=\delta B_{22},$ and $\hat{P}_{i,N+1}:=H_i$. $\hat{\mathbb{E}}_{k}(x_s)$ means the conditional expectation of $x_s(\cdot)$ with respect to $\hat{\mathcal{F}}_k$. Then we have the following result.
		\begin{mylem}
			The discretized form $(x_k, p_{i,k}, q_{i,k})$ of DFBSDEs (\ref{DFBSDEs}) satisfies the following equations:
			\begin{align}
			x_{k+1}&=A_kx_k+B_{11,k}\hat{\mathbb{E}}_{k-d_1-1}(p_{1,k})+B_{12,k}\hat{\mathbb{E}}_{k-d_2-1}(p_{2,k})\nonumber\\ 
            &\quad+\frac{1}{\delta }B_{21,k}\hat{\mathbb{E}}_{k-d_1-1}(\Delta w_kp_{1,k})+\frac{1}{\delta}B_{22,k}\hat{\mathbb{E}}_{k-d_2-1}(\Delta w_kp_{2,k}),\label{x k} \\
			p_{i,k-1}&=\hat{\mathbb{E}}_{k-1}(A'_kp_{i,k})+\hat{Q}_ix_k,\label{p k}\\ 
            p_{i,N}&=H_ix_{N+1}=\hat{P}_{i,N+1}x_{N+1},\quad i=1,2 \label{p N}.
			\end{align}
		\end{mylem}
		
		\noindent\textit{Proof.} Referring to Kailath et al. \cite{KSH2000}, we discretize the backward equation in (\ref{DFBSDEs}) as
		\begin{equation}\label{p k another form}
			p_{i,k}-p_{i,k-1}=-\delta[A'p_{i,k}+\bar{A}'q_{i,k}+Q_ix_k]+\Delta w_kq_{i,k}.
		\end{equation}
		Then multiplying $\Delta w_k$ and taking the conditional expectation with regard to $\hat{\mathcal{F}}_{k-1}$ on both sides, we get
		\begin{align*}
			\hat{\mathbb{E}}_{k-1}[\Delta w_k(p_{i,k}-p_{i,k-1})]=-\hat{\mathbb{E}}_{k-1}[\delta\Delta w_k(A'p_{i,k}+\bar{A}'q_{i,k}+Q_ix_k)]+\hat{\mathbb{E}}_{k-1}(\Delta w_k^2q_{i,k}),
		\end{align*}
		which is 
		\begin{align*}
			\hat{\mathbb{E}}_{k-1}(\Delta w_kp_{i,k})=-\hat{\mathbb{E}}_{k-1}[\delta\Delta w_k(A'p_{i,k}+\bar{A}'q_{i,k}+Q_ix_k)]+\delta q_{i,k}.
		\end{align*}
		Noting that the first term of the right hand of the equality is higher order of $\delta$ due to the fact that $\mathbb{E}(\Delta w_k^2)=\delta$, it will tend to zero with $\delta$ tending to zero at last. The error is allowable so we omit it here, which gives $q_{i,k}=\frac{1}{\delta}\hat{\mathbb{E}}_{k-1}(\Delta w_kp_{i,k})$.
		Then inserting $q_{i,k}=\frac{1}{\delta}\hat{\mathbb{E}}_{k-1}(\Delta w_kp_{i,k})$ into (\ref{p k another form}), and taking the conditional expectation, we get
		\begin{align*}
			p_{i,k-1}=\hat{\mathbb{E}}_{k-1}[(I+\delta A'+\Delta w_k\bar{A}')p_{i,k}]+\delta Q_ix_k=\hat{\mathbb{E}}_{k-1}(A'_kp_{i,k})+\hat{Q}_ix_k,
		\end{align*}
		which is (\ref{p k}). In the same way, we discretize the forward equation in (\ref{DFBSDEs}) as
		\begin{align}
			x_{k+1}-x_k&=(\delta A+\Delta w_k\bar{A})x_k+(\delta B_{11}+\Delta w_k\bar{B}_{11})\hat{\mathbb{E}}_{k-d_1-1}(p_{1,k})+(\delta B_{12}\nonumber\\
            &\quad +\Delta w_k\bar{B}_{12})\hat{\mathbb{E}}_{k-d_2-1}(p_{2,k})+(\delta B_{21}+\Delta w_k\bar{B}_{21})\hat{\mathbb{E}}_{k-d_1-1}(q_{1,k})\nonumber\\ 
            &\quad +(\delta B_{22}+\Delta w_k\bar{B}_{22})\hat{\mathbb{E}}_{k-d_2-1}(q_{2,k}).\label{x k another form}
		\end{align}
		After putting $q_{i,k}=\frac{1}{\delta}\hat{\mathbb{E}}_{k-1}(\Delta w_kp_{i,k})$ into (\ref{x k another form}), we could get (\ref{x k}). $\hfill\qedsymbol$
		
		Next, through the following lemma, we obtain the discrete-time form of explicit solution for (\ref{x k})-(\ref{p N}).
		\begin{mylem}
			If Riccati equations (\ref{Riccati 1})-(\ref{Riccati 5}) admit solutions such that matrices $\hat{\Gamma}_k, \Gamma_k^i, \check{\Gamma}_k, i=1,2$ are invertible, then the solution of DFBSDEs (\ref{x k})-(\ref{p N}) for $k\geq d_1$ are presented as follows:
			\begin{align}
				x_k&=A_{k-1}x_{k-1}+M_{k-1}\hat{\mathbb{E}}_{k-d_1-2}(x_{k-1})+\sum_{m=1}^{d_1-d_2-1}M_{k-1}^m\hat{\mathbb{E}}_{k-d_1-2-m}(x_{k-1})\nonumber\\&\quad+H_{k-1}\hat{\mathbb{E}}_{k-d_2-2}(x_{k-1}),\label{x k explicit}\\
				p_{i,k-1}&=\hat{P}_{i,k}x_k+\sum_{j=0}^{d_1}\hat{P}_{i,k}^{k+j}\hat{\mathbb{E}}_{k-d_1-1+j}(x_k)+\sum_{j=0}^{d_2}\check{P}_{i,k}^{k+j}\hat{\mathbb{E}}_{k-d_2-1+j}(x_k),\label{p k explicit}
			\end{align}
			and $\hat{P}_{i,k}, \hat{P}_{i,k}^{k+m}, \check{P}_{i,k}^{k+m}$ satisfy
			\begin{align}
				\hat{P}_{i,k}&=\hat{A}'\hat{P}_{i,k+1}\hat{A}+\delta\bar{A}'\hat{P}_{i,k+1}\bar{A}+\hat{A}'(\hat{P}_{i,k+1}^{k+1+d_1}+\check{P}_{i,k+1}^{k+1+d_2})\hat{A}+\hat{Q}_i,\label{Riccati 1}\\
				\hat{P}_{i,k}^{k}&=\mathbb{E}\big[(\hat A'\hat{S}_{i,k+1}+\Delta w_k\bar{A}'\hat{P}_{i,k+1})M_k\big],\label{Riccati 2}\\
				\hat{P}_{i,k}^{k+m}&=
				\begin{cases}
					\mathbb{E}\big[(\hat A'\hat{S}^m_{i,k+1}+\Delta w_k\bar{A}'\hat{P}_{i,k+1})M_k^m+\hat{A}'\hat{P}_{i,k+1}^{k+m}(\hat{A}+H_k+\sum_{j=m}^{d_1-d_2-1}M_k^j)\big],\\
                    \hfill 0<m<d_1-d_2,\\
                    \hat{A}'\hat{P}_{i,k+1}^{k+m}\hat{A}, \hfill d_1-d_2 \leq m,
				\end{cases}\label{Riccati 3}\\
				\check{P}_{i,k}^k&=\mathbb{E}\big[(\hat{A}'\check{S}_{i,k+1}+\Delta w_k\bar{A}'\hat{P}_{i,k+1})H_k\big],\label{Riccati 4}\\ 
                \check{P}_{i,k}^{k+m}&=\hat{A}'\check{P}_{i,k+1}^{k+m}\hat{A}.\label{Riccati 5}
			\end{align}	
			In the above, the coefficient matrices are defined as follows:		
			\begin{align}
				M_{k-1}&:=\Bigg\{\bigg[B_{11,k-1}\hat{S}_{1,k}\quad \frac{1}{\delta}B_{21,k-1}\hat{P}_{1,k}\bigg] +\bigg[B_{12,k-1}\check{S}_{2,k}\quad \frac{1}{\delta}B_{22,k-1}\hat{P}_{2,k}\bigg]\check{\Gamma}_k^{-1}G_k\nonumber\\
                &\qquad +\bigg\{\big[B_{12,k-1}\quad 0\big]+\bigg[B_{12,k-1}\check{S}_{2,k}\quad \frac{1}{\delta}B_{22,k-1}\hat{P}_{2,k}\bigg]\check{\Gamma}_k^{-1}
                \begin{bmatrix}
					\hat{B}_{12} & 0 \\	\delta\bar{B}_{12} & 0
				\end{bmatrix}\bigg\}\nonumber\\
                &\qquad \times\bigg[\sum_{j=1}^{d_1-d_2-1}Z_k(j)\hat{P}_{2,k}^{k+j-1}(\Gamma_k^j)^{-1}\bigg]G_k\Bigg\}\hat{\Gamma}_k^{-1}
                \left[\begin{array}{cc}
				      \hat{A} \\ \delta\bar{A}  
				\end{array}\right],\label{M k}\\
				M_{k-1}^m&:=\Bigg\{[B_{12,k-1}\quad 0]+\bigg[B_{12,k-1}\check{S}_{2,k}\quad \frac{1}{\delta}B_{22,k-1}\hat{P}_{2,k}\bigg]\check{\Gamma}_k^{-1}
                \begin{bmatrix}
				      \hat{B}_{12} & 0 \\
				      \delta\bar{B}_{12} & 0
				\end{bmatrix}\Bigg\}Z_k(m)\nonumber\\
                &\quad \times\hat{P}_{2,k}^{k+m-1}(\Gamma_k^m)^{-1}
                \left[\begin{array}{cc}
				      \hat{A} \\ \delta\bar{A}  
				\end{array}\right],\label{M m}\\
				H_{k-1}&:=\bigg[B_{12,k-1}\check{S}_{2,k}\quad\dfrac{1}{\delta}B_{22,k-1}\hat{P}_{2,k}\bigg]\check{\Gamma}_k^{-1}
                \left[\begin{array}{cc}
				      \hat{A} \\ \delta\bar{A}  
				\end{array}\right],\label{H k}\\
				\hat{\Gamma}_k&:=\left[\begin{array}{cc}
				     I-\hat{B}_{11}\hat{S}_{1,k}-\hat{B}_{12}\hat{S}_{2,k} & -\frac{1}{\delta}(\hat{B}_{21}\hat{P}_{1,k}+\hat{B}_{22}\hat{P}_{2,k}) \\ 
                     -\delta(\bar{B}_{11}\hat{S}_{1,k}+\bar{B}_{12}\hat{S}_{2,k}) & I-\bar{B}_{21}\hat{P}_{1,k}-\bar{B}_{22}\hat{P}_{2,k} 
				\end{array}\right],\label{Gamma k hat}
                \end{align}
				\begin{align}
				\Gamma_k^m&=\left[\begin{array}{cc}
						I-\hat{B}_{12}S_{2,k}^m & -\frac{1}{\delta}\hat{B}_{22}\hat{P}_{2,k} \\ -\delta\bar{B}_{12}S_{2,k}^m & I-\bar{B}_{22}\hat{P}_{2,k} 
					\end{array}\right],\quad
				\check{\Gamma}_k=\left[\begin{array}{cc}
						I-\hat{B}_{12}\check{S}_{2,k} & -\frac{1}{\delta}\hat{B}_{22}\hat{P}_{2,k} \\ -\delta\bar{B}_{12}\check{S}_{2,k} & I-\bar{B}_{22}\hat{P}_{2,k} 
					\end{array}\right],\label{Gamma k m}\\
				G_k&:=\left[\begin{array}{cc}
						\hat{B}_{11}\hat{S}_{1,k} & \frac{1}{\delta}\hat{B}_{21}\hat{P}_{1,k} \\ \delta\bar{B}_{11}\hat{S}_{1,k} & \bar{B}_{21}\hat{P}_{1,k}
						\end{array}\right],\quad
					\hat{S}_{i,k}:=\hat{P}_{i,k}+\sum_{j=0}^{d_1}\hat{P}_{i,k}^{k+j}+\sum_{j=0}^{d_2}\check{P}_{i,k}^{k+j},\label{G k}\\
				S_{i,k}^m&:=\hat{P}_{i,k}+\sum_{j=m-1}^{d_1}\hat{P}_{i,k}^{k+j}+\sum_{j=0}^{d_2}\check{P}_{i,k}^{k+j},\quad
				\check{S}_{i,k}:=\hat{P}_{i,k}+\sum_{j=d_1-d_2-1}^{d_1}\hat{P}_{i,k}^{k+j}+\sum_{j=0}^{d_2}\check{P}_{i,k}^{k+j},\label{S k m hat}
			\end{align}
			with $p_{i,N}=\hat{P}_{i,N+1}x_{N+1}, \hat{P}_{i,k}^m=\check{P}_{i,k}^m=0$ if $m>N$, $i=1,2$, $M_k^m=0$, if $m\geq d_1-d_2$ or $k+m>N$.
		\end{mylem}
		\noindent\textit{Proof}. Considering $k=N$, we combine (\ref{x k}) and (\ref{p N}),
		\begin{align}\label{x N}
			x_{N+1}=&A_Nx_N+B_{11,N}\hat{\mathbb{E}}_{N-d_1-1}(p_{1,N})+B_{12,N}\hat{\mathbb{E}}_{N-d_2-1}(p_{2,N})\nonumber\\ 
                    &+\frac{1}{\delta }B_{21,N}\hat{\mathbb{E}}_{N-d_1-1}(\Delta w_Np_{1,N})+\frac{1}{\delta}B_{22,N}\hat{\mathbb{E}}_{N-d_2-1}(\Delta w_Np_{2,N})\nonumber\\
			       =&A_Nx_N+B_{11,N}\hat{\mathbb{E}}_{N-d_1-1}(\hat{P}_{1,N+1}x_{N+1})+B_{12,N}\hat{\mathbb{E}}_{N-d_2-1}(\hat{P}_{2,N+1}x_{N+1})\nonumber\\ 
                    &+\frac{1}{\delta }B_{21,N}\hat{\mathbb{E}}_{N-d_1-1}(\Delta w_N\hat{P}_{1,N+1}x_{N+1})+\frac{1}{\delta}B_{22,N}\hat{\mathbb{E}}_{N-d_2-1}(\Delta w_N\hat{P}_{2,N+1}x_{N+1}).
		\end{align}
		Taking the conditional expectation with regard to $\hat{\mathcal{F}}_{N-d_1-1}$, it yields
		\begin{align}\label{x N11}
			&\big(I-\hat{B}_{11}\hat{P}_{1,N+1}-\hat{B}_{12}\hat{P}_{2,N+1}\big)\hat{\mathbb{E}}_{N-d_1-1}(x_{N+1})\nonumber\\
			&=\hat{A}\hat{\mathbb{E}}_{N-d_1-1}(x_N)+\frac{1}{\delta}(\hat{B}_{21}\hat{P}_{1,N+1}+\hat{B}_{22}\hat{P}_{2,N+1})\hat{\mathbb{E}}_{N-d_1-1}(\Delta w_Nx_{N+1}).
		\end{align}
		Multiplying $\Delta w_N$ and taking the conditional expectation on both sides, we could get 
		\begin{align}\label{x N12}
			&\big(I-\bar{B}_{21}\hat{P}_{1,N+1}-\bar{B}_{22}\hat{P}_{2,N+1}\big)\hat{\mathbb{E}}_{N-d_1-1}(\Delta w_Nx_{N+1})\nonumber\\
			&=\delta\bar{A}\hat{\mathbb{E}}_{N-d_1-1}(x_N)+\delta(\bar{B}_{11}\hat{P}_{1,N+1}+\bar{B}_{12}\hat{P}_{2,N+1})\hat{\mathbb{E}}_{N-d_1-1}(x_{N+1}).
		\end{align}
		Combining (\ref{x N11})-(\ref{x N12}) and considering the form of $\hat{\Gamma}_{N+1}$, we have 
		\begin{align}\label{x N1}
			\left[\begin{array}{cc}
			\hat{\mathbb{E}}_{N-d_1-1}(x_{N+1})	 \\ \hat{\mathbb{E}}_{N-d_1-1}(\Delta w_Nx_{N+1})
			\end{array}\right]=\hat{\Gamma}_{N+1}^{-1}\left[\begin{array}{cc}
			\hat{A}	 \\ \delta\bar{A}\end{array}\right]\hat{\mathbb{E}}_{N-d_1-1}(x_N).
		\end{align}
		Similarly, take the conditional expectation with regard to $\hat{\mathcal{F}}_{N-d_2-1}$,
		\begin{align}\label{x N21}
			&(I-\hat{B}_{12}\hat{P}_{2,N+1})\hat{\mathbb{E}}_{N-d_2-1}(x_{N+1})=\hat{A}\hat{\mathbb{E}}_{N-d_2-1}(x_N)+\hat{B}_{11}\hat{P}_{1,N+1}\hat{\mathbb{E}}_{N-d_1-1}(x_{N+1})\nonumber\\
            &\qquad +\frac{1}{\delta}\hat{B}_{21}\hat{P}_{1,N+1}\hat{\mathbb{E}}_{N-d_1-1}(\Delta w_Nx_{N+1})+\frac{1}{\delta}\hat{B}_{22}\hat{P}_{2,N+1}\hat{\mathbb{E}}_{N-d_2-1}(\Delta w_Nx_{N+1}).
		\end{align}
		Multiplying $\Delta w_N$ and taking the conditional expectation, we have 
		\begin{align}\label{x N22}
			&(I-\bar{B}_{22}\hat{P}_{2,N+1})\hat{\mathbb{E}}_{N-d_2-1}(\Delta w_Nx_{N+1})=\delta\bar{A}\hat{\mathbb{E}}_{N-d_2-1}(x_N)+\delta\bar{B}_{11}\hat{P}_{1,N+1}\hat{\mathbb{E}}_{N-d_1-1}(x_{N+1})\nonumber\\
            &\qquad +\bar{B}_{21}\hat{P}_{1,N+1}\hat{\mathbb{E}}_{N-d_1-1}(\Delta w_Nx_{N+1})+\delta\bar{B}_{12}\hat{P}_{2,N+1}\hat{\mathbb{E}}_{N-d_2-1}(x_{N+1}).
		\end{align}
		Combining (\ref{x N21})-(\ref{x N22}) and considering the form of $\check{\Gamma}_{N+1}$, it generates 
		\begin{align}\label{x N2}
			&\left[\begin{array}{cc}
				\hat{\mathbb{E}}_{N-d_2-1}(x_{N+1})	 \\ \hat{\mathbb{E}}_{N-d_2-1}(\Delta w_Nx_{N+1})
			\end{array}\right]\nonumber\\ 
            &=\check{\Gamma}_{N+1}^{-1}\left[\begin{array}{cc}
				\hat{A}	 \\ \delta\bar{A}\end{array}\right]\hat{\mathbb{E}}_{N-d_2-1}(x_N)
             +\check{\Gamma}_{N+1}^{-1}G_{N+1}\left[\begin{array}{cc}
				\hat{\mathbb{E}}_{N-d_1-1}(x_{N+1})	 \\ \hat{\mathbb{E}}_{N-d_1-1}(\Delta w_Nx_{N+1})
				\end{array}\right].
		\end{align}
		Introducing (\ref{x N1}) and (\ref{x N2}) into (\ref{x N}), the following equation could be obtained immediately,  
		\begin{align}\label{x N3}
			x_{N+1}=A_Nx_N+M_N\hat{\mathbb{E}}_{N-d_1-1}(x_N)+H_N\hat{\mathbb{E}}_{N-d_2-1}(x_N).
		\end{align}
		Then combining (\ref{p k}), (\ref{p N}) and (\ref{x N3}), we have 
		\begin{align*}
			p_{i,N-1}&=\hat{\mathbb{E}}_{N-1}(A'_Np_{i,N})+\hat{Q}_ix_N\nonumber\\
                     &=\hat{\mathbb{E}}_{N-1}\{A'_N\hat{P}_{i,N+1}[A_Nx_N+M_N\hat{\mathbb{E}}_{N-d_1-1}(x_N)+H_N\hat{\mathbb{E}}_{N-d_2-1}(x_N)]+\hat{Q}_ix_N\nonumber\\
                     &=\hat{P}_{i,N}x_N+\hat{P}_{i,N}^N\hat{\mathbb{E}}_{N-d_1-1}(x_N)+\check{P}_{i,N}^N\hat{\mathbb{E}}_{N-d_2-1}(x_N).
		\end{align*}
		Thus (\ref{x k explicit})-(\ref{p k explicit}) are true for $k=N$.
		Then taking any $s$ with $d_1<s<N$. Assuming (\ref{x k explicit})-(\ref{p k explicit}) hold at $k\geq s+1$, we will prove that they also hold for $k=s$. Combining (\ref{p k}), (\ref{x k explicit}) and (\ref{p k explicit}), we have
		\begin{align*}
			p_{i,s-1}&=\hat{\mathbb{E}}_{s-1}(A'_sp_{i,s})+\hat{Q}_ix_s\\
                     &=\hat{\mathbb{E}}_{s-1}\bigg\{A'_s\bigg[\hat{P}_{i,s+1}x_{s+1}+\sum_{j=0}^{d_1}\hat{P}_{i,s+1}^{s+1+j}\hat{\mathbb{E}}_{s-d_1+j}(x_{s+1})+\sum_{j=0}^{d_2}\check{P}_{i,s+1}^{s+1+j}\hat{\mathbb{E}}_{s-d_2+j}(x_{s+1})\bigg]\bigg\}
                      +\hat{Q}_ix_s.
		\end{align*}
		Noting that for any $s>N, \hat{P}_{i,k}^s=\check{P}_{i,k}^s=0$, $i=1,2$, we arrive at the below equality:
		\begin{align*}
			p_{i,s-1}&=\big[\hat{A}'\hat{P}_{i,s+1}\hat{A}+\delta\bar{A}'\hat{P}_{i,s+1}\bar{A}+\hat{Q}_i+\hat{A}'(\hat{P}_{i,s+1}^{s+d_1+1}+\check{P}_{i,s+1}^{s+d_2+1})\hat{A}\big]x_s\\
                     &\quad +\mathbb{E}\bigg\{\bigg[A'_s\hat{P}_{i,s+1}+\hat{A}'\bigg(\sum_{j=0}^{d_1}\hat{P}_{i,s+1}^{s+1+j}+\sum_{j=0}^{d_2}\check{P}_{i,s+1}^{s+1+j}\bigg)\bigg]M_s\bigg\}\hat{\mathbb{E}}_{s-d_1-1}(x_s)\\
                     &\quad +\mathbb{E}\bigg\{\bigg[A'_s\hat{P}_{i,s+1}+\hat{A}'\bigg(\sum_{j=1}^{d_1}\hat{P}_{i,s+1}^{s+1+j}+\sum_{j=0}^{d_2}\check{P}_{i,s+1}^{s+1+j}\bigg)\bigg]M_s^1\\
                     &\qquad\qquad +\hat{A}'\hat{P}_{i,s+1}^{s+1}\bigg(\hat{A}+\sum_{m=1}^{d_1-d_2-1}M_s^m+H_s\bigg)\bigg\}\hat{\mathbb{E}}_{s-d_1}(x_s)+\dots \\
                     &\quad +\mathbb{E}\bigg\{\bigg[A'_s\hat{P}_{i,s+1}+\hat{A}'\bigg(\sum_{j=d_1-d_2-1}^{d_1}\hat{P}_{i,s+1}^{s+1+j}+\sum_{j=0}^{d_2}\check{P}_{i,s+1}^{s+1+j}\bigg)\bigg]M_s^{d_1-d_2-1}\\
                     &\qquad\qquad +\hat{A}'\hat{P}_{i,s+1}^{s+d_1-d_2-1}\big(\hat{A}+M_s^{d_1-d_2-1}+H_s\big)\bigg\}\hat{\mathbb{E}}_{s-d_2-2}(x_s)\\
        \end{align*}
        \begin{align*}
                     &\quad +\mathbb{E}\bigg\{\bigg[A'_s\hat{P}_{i,s+1}+\hat{A}'\bigg(\sum_{j=d_1-d_2-1}^{d_1}\hat{P}_{i,s+1}^{s+1+j}+\sum_{j=0}^{d_2}\check{P}_{i,s+1}^{s+1+j}\bigg)\bigg]H_s\\
                     &\qquad\qquad +\hat{A}'\hat{P}_{i,s+1}^{s+d_1-d_2}\hat{A}\bigg\}\hat{\mathbb{E}}_{s-d_2-1}(x_s)+\big(\hat{A}'\hat{P}_{i,s+1}^{s+d_1-d_2+1}\hat{A}\big)\hat{\mathbb{E}}_{s-d_2}(x_s)\\
                     &\quad +(\hat{A}'\check{P}_{i,s+1}^{s+1}\hat{A})\hat{\mathbb{E}}_{s-d_2}(x_s)+\dots+(\hat{A}'\hat{P}_{i,s+1}^{s+d_1}\hat{A})\hat{\mathbb{E}}_{s-1}(x_s)+(\hat{A}'\check{P}_{i,s+1}^{s+d_2}\hat{A})\hat{\mathbb{E}}_{s-1}(x_s)\\
                     &=\hat{P}_{i,s}x_s+\sum_{j=0}^{d_1}\hat{P}_{i,s}^{s+j}\hat{\mathbb{E}}_{s-d_1-1+j}(x_s)+\sum_{j=0}^{d_2}\check{P}_{i,s}^{s+j}\hat{\mathbb{E}}_{s-d_2-1+j}(x_s).
		\end{align*}
		Using (\ref{x k}) and (\ref{p k explicit}), and noting the form of $\hat{S}_{i,k}, S_{i,k}^m, \check{S}_{i,k}$, we have
		\begin{align}\label{x s}
			x_{s}=&A_{s-1}x_{s-1}+B_{11,s-1}\hat{\mathbb{E}}_{s-d_1-2}(p_{1,s-1})+B_{12,s-1}\hat{\mathbb{E}}_{s-d_2-2}(p_{2,s-1})\nonumber\\ 
			      &+\frac{1}{\delta }B_{21,s-1}\hat{\mathbb{E}}_{s-d_1-2}(\Delta w_{s-1}p_{1,s-1})+\frac{1}{\delta}B_{22,s-1}\hat{\mathbb{E}}_{s-d_2-2}(\Delta w_{s-1}p_{2,s-1})\nonumber\\
                 =&A_{s-1}x_{s-1}+B_{11,s-1}\hat{S}_{1,s}\hat{\mathbb{E}}_{s-d_1-2}(x_s)+\sum_{j=0}^{d_1-d_2-2}B_{12,s-1}\hat{P}_{2,s}^{s+j}\hat{\mathbb{E}}_{s-d_1-1+j}(x_s)\nonumber\\
			      &+B_{12,s-1}\check{S}_{2,s}\hat{\mathbb{E}}_{s-d_2-2}(x_s)+\frac{1}{\delta }B_{21,s-1}\hat{P}_{1,s}\hat{\mathbb{E}}_{s-d_1-2}(\Delta w_{s-1}x_s)\nonumber\\
			      &+\frac{1}{\delta }B_{22,s-1}\hat{P}_{2,s}\hat{\mathbb{E}}_{s-d_2-2}(\Delta w_{s-1}x_s).
		\end{align}
		For the sake of simplicity, the subsequent discussion will take $d_1-d_2=3$ as an illustration. Just like the derivation processes of (\ref{x N1}) and (\ref{x N2}), we can obtain 
		\begin{align}\label{x s1}
		&\left[\begin{array}{cc}
			I-\hat{B}_{11}\hat{S}_{1,s}-\hat{B}_{12}\hat{S}_{2,s} & -\dfrac{1}{\delta}(\hat{B}_{21}\hat{P}_{1,s}+\hat{B}_{22}\hat{P}_{2,s})	 \\ 
            -\delta(\bar{B}_{11}\hat{S}_{1,s}+\bar{B}_{12}\hat{S}_{2,s}) & I-(\bar{B}_{21}\hat{P}_{1,s}+\bar{B}_{22}\hat{P}_{2,s})
		\end{array}\right]
        \left[\begin{array}{cc}
			\hat{\mathbb{E}}_{s-d_1-2}(x_s)	 \\ \hat{\mathbb{E}}_{s-d_1-2}(\Delta w_{s-1}x_s)
		\end{array}\right]\nonumber\\
		&=\left[\begin{array}{cc}
		\hat{A}	 \\ 	\delta\bar{A}\end{array}\right]\hat{\mathbb{E}}_{s-d_1-2}(x_{s-1}),
		\end{align}
		\begin{align}\label{x s2}
		&\left[\begin{array}{cc}
			I-\hat{B}_{12}S_{2,s}^1 & -\dfrac{1}{\delta}\hat{B}_{22}\hat{P}_{2,s} \\ -\delta\bar{B}_{12}S_{2,s}^1 & I-\bar{B}_{22}\hat{P}_{2,s}
		\end{array}\right]\left[\begin{array}{cc}
			\hat{\mathbb{E}}_{s-d_1-1}(x_s)	 \\ \hat{\mathbb{E}}_{s-d_1-1}(\Delta w_{s-1}x_s)
		\end{array}\right]\nonumber\\
		&=\left[\begin{array}{cc}
			\hat{A}	 \\ 	\delta\bar{A}\end{array}\right]\hat{\mathbb{E}}_{s-d_1-1}(x_{s-1})+G_s\left[\begin{array}{cc}
			\hat{\mathbb{E}}_{s-d_1-2}(x_s)	 \\ \hat{\mathbb{E}}_{s-d_1-2}(\Delta w_{s-1}x_s)
			\end{array}\right],
		\end{align}		
		\begin{align}\label{x s3}
		&\left[\begin{array}{cc}
			I-\hat{B}_{12}S_{2,s}^2 & -\dfrac{1}{\delta}\hat{B}_{22}\hat{P}_{2,s} \\ -\delta\bar{B}_{12}S_{2,s}^2 & I-\bar{B}_{22}\hat{P}_{2,s}
		\end{array}\right]\left[\begin{array}{cc}
			\hat{\mathbb{E}}_{s-d_1}(x_s)	 \\ \hat{\mathbb{E}}_{s-d_1}(\Delta w_{s-1}x_s)
		\end{array}\right]\nonumber\\
		&=\left[\begin{array}{cc}
				\hat{A}	 \\ 	\delta\bar{A}\end{array}\right]\hat{\mathbb{E}}_{s-d_1}(x_{s-1})+G_s\left[\begin{array}{cc}
				\hat{\mathbb{E}}_{s-d_1-2}(x_s)	 \\ \hat{\mathbb{E}}_{s-d_1-2}(\Delta w_{s-1}x_s)
			\end{array}\right]\nonumber\\
		&\quad+\left[\begin{array}{cc}
				\hat{B}_{12} & 0 \\ 
				\delta\bar{B}_{12} & 0 \end{array}\right]\hat{P}_{2,s}^s
			\left[\begin{array}{cc}
				\hat{\mathbb{E}}_{s-d_1-1}(x_s)	 \\ \hat{\mathbb{E}}_{s-d_1-1}(\Delta w_{s-1}x_s)
			\end{array}\right],
		\end{align}
		\begin{align}\label{x s4}
		&\left[\begin{array}{cc}
				I-\hat{B}_{12}\check{S}_{2,s} & -\dfrac{1}{\delta}\hat{B}_{22}\hat{P}_{2,s} \\ -\delta\bar{B}_{12}\check{S}_{2,s} & I-\bar{B}_{22}\hat{P}_{2,s}
			\end{array}\right]\left[\begin{array}{cc}
				\hat{\mathbb{E}}_{s-d_2-2}(x_s)	 \\ \hat{\mathbb{E}}_{s-d_2-2}(\Delta w_{s-1}x_s)
			\end{array}\right]\nonumber\\
		&=\left[\begin{array}{cc}
				\hat{A}	 \\ \delta\bar{A}\end{array}
			\right]\hat{\mathbb{E}}_{s-d_2-2}(x_{s-1})+G_s\left[\begin{array}{cc}
				\hat{\mathbb{E}}_{s-d_1-2}(x_s)	 \\ \hat{\mathbb{E}}_{s-d_1-2}(\Delta w_{s-1}x_s)
			\end{array}\right]\nonumber\\
		&\quad+\left[\begin{array}{cc}
			\hat{B}_{12} & 0 \\ \delta\bar{B}_{12} & 0\end{array}\right]\Bigg(\hat{P}_{2,s}^s
			\left[\begin{array}{cc}
				\hat{\mathbb{E}}_{s-d_1-1}(x_s)	 \\ \hat{\mathbb{E}}_{s-d_1-1}(\Delta w_{s-1}x_s)
			\end{array}\right]+\hat{P}_{2,s}^{s+1}
			\left[\begin{array}{cc}
			\hat{\mathbb{E}}_{s-d_1}(x_s)	 \\ \hat{\mathbb{E}}_{s-d_1}(\Delta w_{s-1}x_s)
			\end{array}\right]\Bigg).
		\end{align}
		Putting (\ref{x s1}), (\ref{x s2}), (\ref{x s3}) and (\ref{x s4}) into (\ref{x s}), $x_s$ could be rewritten as follows:
		\begin{align*}
			x_s&=A_{s-1}x_{s-1}+\left[\begin{array}{cc}
				B_{11,s-1}\hat{S}_{1,s} & \frac{1}{\delta}B_{21,s-1}\hat{P}_{1,s}\end{array}\right]
				\left[\begin{array}{cc}
						\hat{\mathbb{E}}_{s-d_1-2}(x_s)	 \\ \hat{\mathbb{E}}_{s-d_1-2}(\Delta w_{s-1}x_s)
					\end{array}\right]\\
			&\quad+\left[\begin{array}{cc}
				B_{12,s-1} & 0
				\end{array}\right]\Bigg(\hat{P}_{2,s}^{s}\left[\begin{array}{cc}
					\hat{\mathbb{E}}_{s-d_1-1}(x_s)	 \\ \hat{\mathbb{E}}_{s-d_1-1}(\Delta w_{s-1}x_s)
				\end{array}\right]+\hat{P}_{2,s}^{s+1}\left[\begin{array}{cc}
					\hat{\mathbb{E}}_{s-d_1}(x_s)	 \\ \hat{\mathbb{E}}_{s-d_1}(\Delta w_{s-1}x_s)
				\end{array}\right]\Bigg)\\
			&\quad+\left[\begin{array}{cc}
						B_{12,s-1}\check{S}_{2,s} & \frac{1}{\delta}B_{22,s-1}\hat{P}_{2,s}
					\end{array}\right]
					\left[\begin{array}{cc}
						\hat{\mathbb{E}}_{s-d_2-2}(x_s)	 \\ \hat{\mathbb{E}}_{s-d_2-2}(\Delta w_{s-1}x_s)
					\end{array}\right]\\
			&=A_{s-1}x_{s-1}
			   		+\Bigg\{\bigg[B_{11,s-1}\hat{S}_{1,s}\quad\frac{1}{\delta}B_{21,s-1}\hat{P}_{1,s}\bigg]+\bigg[B_{12,s-1}\check{S}_{2,s}\quad \frac{1}{\delta}B_{22,s-1}\hat{P}_{2,s}\bigg]\check{\Gamma}_s^{-1}G_s\\
			&\qquad+\bigg\{\big[B_{12,s-1}\quad 	0\big]+\bigg[B_{12,s-1}\check{S}_{2,s}\quad \frac{1}{\delta}B_{22,s-1}\hat{P}_{2,s}\bigg]\check{\Gamma}_s^{-1}\begin{bmatrix}
			   		\hat{B}_{12} & 0 \\	\delta\bar{B}_{12} & 0
			    \end{bmatrix}\bigg\}\bigg[\hat{P}_{2,s}^{s+1}(\Gamma_s^2)^{-1}\\
			&\qquad+\Big(I+\hat{P}_{2,s}^{s+1}(\Gamma_s^2)^{-1}\begin{bmatrix}
				\hat{B}_{12} & 0 \\	\delta\bar{B}_{12} & 0
			\end{bmatrix}\Big)\hat{P}_{2,s}^s(\Gamma_s^1)^{-1}\bigg]G_s\Bigg\}\hat{\Gamma}_s^{-1}
			   \left[\begin{array}{cc}
			   	 \hat{A} \\ \delta\bar{A}  
			   \end{array}\right]\hat{\mathbb{E}}_{s-d_1-2}(x_{s-1})\\
			&\quad +\Bigg\{[B_{12,s-1}\quad 0]+\bigg[B_{12,s-1}\check{S}_{2,s}\quad \frac{1}{\delta}B_{22,s-1}\hat{P}_{2,s}\bigg]\check{\Gamma}_s^{-1}\begin{bmatrix}
				\hat{B}_{12} & 0 \\
				\delta\bar{B}_{12} & 0
			\end{bmatrix}\Bigg\}\\
			&\quad \times\Big(I+\hat{P}_{2,s}^{s+1}(\Gamma_s^2)^{-1}\begin{bmatrix}
				\hat{B}_{12} & 0 \\
				\delta\bar{B}_{12} & 0
			\end{bmatrix}\Big)\hat{P}_{2,s}^s(\Gamma_s^1)^{-1}
			\left[\begin{array}{cc}
				\hat{A} \\ \delta\bar{A}  
			\end{array}\right]\hat{\mathbb{E}}_{s-d_1-1}(x_{s-1})\\
			&\quad+\Bigg\{[B_{12,s-1}\quad 0]+\bigg[B_{12,s-1}\check{S}_{2,s}\quad \frac{1}{\delta}B_{22,s-1}\hat{P}_{2,s}\bigg]\check{\Gamma}_s^{-1}\begin{bmatrix}
				\hat{B}_{12} & 0 \\
				\delta\bar{B}_{12} & 0
			\end{bmatrix}\Bigg\}\hat{P}_{2,s}^{s+1}(\Gamma_s^2)^{-1}\\
			&\quad \times
				\left[\begin{array}{cc}
					\hat{A} \\ \delta\bar{A}  
				\end{array}\right]\hat{\mathbb{E}}_{s-d_1}(x_{s-1})+\bigg[B_{12,s-1}\check{S}_{2,s}\quad\dfrac{1}{\delta}B_{	22,s-1}\hat{P}_{2,s}\bigg]\check{\Gamma}_s^{-1}
				\left[\begin{array}{cc}
					\hat{A} \\ \delta\bar{A}  
				\end{array}\right]\hat{\mathbb{E}}_{s-d_2-2}(x_{s-1}).
		\end{align*}
		Denoting 
		\begin{align}\label{Z k}
			Z_s(1):=I,\quad Z_s(2):=I+\hat{P}_{2,s}^{s+1}(\Gamma_s^2)^{-1}\begin{bmatrix}
				\hat{B}_{12} & 0 \\
				\delta\bar{B}_{12} & 0
			\end{bmatrix},
		\end{align}
		we could generate
		\begin{align}
			x_s&=A_{s-1}x_{s-1}+M_{s-1}\hat{\mathbb{E}}_{s-d_1-2}(x_{s-1})+M_{s-1}^1\hat{\mathbb{E}}_{s-d_1-1}(x_{s-1})\nonumber\\
			   &\quad +M_{s-1}^2\hat{\mathbb{E}}_{s-d_1}(x_{s-1})+H_{s-1}\hat{\mathbb{E}}_{s-d_2-2}(x_{s-1}).\label{x s final}
		\end{align}
		Therefore, (\ref{x k explicit})-(\ref{Riccati 3}) hold for $k=s$.$\hfill\qedsymbol$
		\begin{Remark}
			We would like to elaborate that, because of the difference of two time delays, $Z_k(\cdot)$ appears in $M_{k-1}, M_{k-1}^m$ of (\ref{M k}), (\ref{M m}), the coefficient matrices of $x_k$. Though in the above proof, for $d_1-d_2=3$, $Z_k(\cdot)$ has some representations in (\ref{Z k}), it will become more and more complex as the difference of time delays increases. However, according to its definition, it would tend to identity matrix with $\delta \rightarrow 0$. Thus, we do not pay too much attention to its specific form, in (\ref{M k}), (\ref{M m}).			
		\end{Remark}

		Lemma 3.4 generates the discrete-time form of explicit solution of DFBSDEs (\ref{x k})-(\ref{p N}), which are depending on the solutions of Riccati equations. In the following result, we let $\delta\rightarrow0 $ to obtain the continuous-time form of explicit solution.
		
		\begin{mythm}
			If Riccati equations (\ref{Riccati 1-1})-(\ref{Riccati 5-5}) admit solutions such that matrices $I-\bar{B}_{21}P_1(t)-\bar{B}_{22}P_2(t), I-\bar{B}_{22}P_2(t)$ are invertible for all $t\in[0,T]$, then (\ref{DFBSDEs}) can be uniquely solved, and the explicit solutions satisfy the following equations, for $i=1,2$,
			\begin{align}
				p_i(t)&=P_i(t)x(t)+\int_{0}^{h_1}\hat{P}_i(t,t+\theta)\mathbb{E}_{t-h_1+\theta}[x(t)]d\theta+\int_{0}^{h_2}\check{P}_i(t,t+\theta)\mathbb{E}_{t-h_2+\theta}[x(t)]d\theta,\label{p i}\\
				q_i(t)&=P_i(t)\Big\{\bar{A}x(t)+\bar{B}_{11}\mathbb{E}_{t-h_1}[p_1(t)]+\bar{B}_{12}\mathbb{E}_{t-h_2}[p_2(t)]\nonumber\\
                      &\quad +\bar{B}_{21}\mathbb{E}_{t-h_1}[q_1(t)]+\bar{B}_{22}\mathbb{E}_{t-h_2}[q_2(t)]\Big\},\label{q i}\\
	             dx(t)&=\bigg\{Ax(t)+\Big\{\big[B_{11}+B_{22}P_2(t)(I-\bar{B}_{22}P_2(t))^{-1}\bar{B}_{11}\big]\hat{S}_1(t)\nonumber\\
                      &\qquad +\big[B_{21}+B_{22}P_2(t)(I-\bar{B}_{22}P_2(t))^{-1}\bar{B}_{21}\big]P_1(t)(I-\bar{B}_{21}P_1(t)\nonumber\\
                      &\qquad -\bar{B}_{22}P_2(t))^{-1}\big[\bar{B}_{11}\hat{S}_1(t)+\bar{B}_{12}\hat{S}_2(t)+\bar{A}\big]\Big\}\mathbb{E}_{t-h_1}[x(t)]\nonumber\\
                      &\qquad +\int_0^{h_1-h_2}\big[B_{12}+B_{22}P_2(t)(I-\bar{B}_{22}P_2(t))^{-1}\bar{B}_{12}\big]\hat{P}_2(t,t+\theta)\mathbb{E}_{t-h_1+\theta}[x(t)]d\theta\nonumber\\ 
                      &\qquad +\big[B_{12}\check{S}_2(t)+B_{22}P_2{t}(I-\bar{B}_{22}P_2(t))^{-1}(\bar{B}_{12}\check{S}_2(t)+\bar{A})\big]\mathbb{E}_{t-h_2}[x(t)]\bigg\}dt\nonumber\\
			          &\quad +\bigg\{\bar{A}x(t)+\Big\{\big[\bar{B}_{11}+\bar{B}_{22}P_2(t)(I-\bar{B}_{22}P_2(t))^{-1}\bar{B}_{11}\big]\hat{S}_1(t)\nonumber\\
                      &\qquad +\big[\bar{B}_{21}+\bar{B}_{22}P_2(t)(I-\bar{B}_{22}P_2(t))^{-1}\bar{B}_{21}\big]P_1(t)(I-\bar{B}_{21}P_1(t)\nonumber\\
                      &\qquad -\bar{B}_{22}P_2(t))^{-1}\big[\bar{B}_{11}\hat{S}_1(t)+\bar{B}_{12}\hat{S}_2(t)+\bar{A}\big]\Big\}\mathbb{E}_{t-h_1}[x(t)]\nonumber\\
                      &\qquad +\int_0^{h_1-h_2}\big[\bar{B}_{12}+\bar{B}_{22}P_2(t)(I-\bar{B}_{22}P_2(t))^{-1}\bar{B}_{12}\big]\hat{P}_2(t,t+\theta)\mathbb{E}_{t-h_1+\theta}[x(t)]d\theta\nonumber\\
                      &\qquad +\big[\bar{B}_{12}\check{S}_2(t)+\bar{B}_{22}P_2{t}(I-\bar{B}_{22}P_2(t))^{-1}(\bar{B}_{12}\check{S}_2(t)+\bar{A})\big]\mathbb{E}_{t-h_2}[x(t)]\Big\}dw(t),\label{x explicit}
			\end{align}
			where $\hat{P}_i(t), \hat{P}_i(t,s)$ and $\check{P}_i(t,s)$ satisfying the following equalities:
			\begin{align}
				 -\dot{P}_i(t)&=A'P_i(t)+P_i(t)A+\bar{A}'P_i(t)\bar{A}+Q_i+\hat{P}_i(t,t+h_1)+\check{P}_i(t,t+h_2),\label{Riccati 1-1}\\
				\hat{P}_i(t,t)&=\big[\hat{S}_i(t)B_{11}+\bar{A}'P_i(t)\bar{B}_{11}\big]\hat{S}_1(t)+\big[\hat{S}_i(t)B_{21}+\bar{A}'P_i(t)\bar{B}_{21}\big]P_1(t)\nonumber\\
				              &\quad \times(I-\bar{B}_{21}P_1(t)-\bar{B}_{22}P_2(t))^{-1}\big[\bar{B}_{11}\hat{S}_1(t)+\bar{B}_{12}\hat{S}_2(t)+\bar{A}\big]\nonumber\\
				              &\quad +\big[\hat{S}_i(t)B_{22}+\bar{A}'P_i(t)\bar{B}_{22}\big]P_2(t)(I-\bar{B}_{22}P_2(t))^{-1}\Big\{\bar{B}_{11}\hat{S}_1(t)\nonumber\\
				              &\quad +\bar{B}_{21}P_1(t)(I-\bar{B}_{21}P_1(t)-\bar{B}_{22}P_2(t))^{-1}\big[\bar{B}_{11}\hat{S}_1(t)+\bar{B}_{12}\hat{S}_2(t)+\bar{A}\big]\Big\},\label{Riccati 2-2}\\
				\check{P}_i(t,t)&=\big[\check{S}_i(t)B_{12}+\bar{A}'P_i(t)\bar{B}_{12}\big]\check{S}_2(t)+\big[\check{S}_i(t)B_{22}+\bar{A}'P_i(t)\bar{B}_{22}\big]P_2(t)\nonumber\\
				              &\quad \times(I-\bar{B}_{22}P_2(t))^{-1}\big[\bar{B}_{12}\check{S}_2(t)+\bar{A}\big],\label{Riccati 3-3}\\
				-\partial t\hat{P}_i(t,s)&=\begin{cases}
					\hat{S}_i(t)\big[B_{12}+B_{22}P_2(t)(I-\bar{B}_{22}P_2(t))^{-1}\bar{B}_{12}\big]\hat{P}_2(t,s)\\
					\quad +\bar{A}'P_i(t)\big[\bar{B}_{12}+\bar{B}_{22}P_2(t)(I-\bar{B}_{22}P_2(t))^{-1}\bar{B}_{12}\big]\hat{P}_2(t,s)\\
                    \quad +\hat{P}_i(t,s)\big[B_{12}\check{S}_2(t)+B_{22}\hat{P}_2(t)(I-\bar{B}_{22}P_2(t))^{-1}\\
                    \quad \times(\bar{B}_{12}\check{S}_2(t)+\bar{A})\big]+A'\hat{P}_i(t,s)+\hat{P}_i(t,s)A,\quad 0<s-t<h_1-h_2,\\
					A'\hat{P}_i(t,s)+\hat{P}_i(t,s)A,\hfill h_1-h_2 \leq s-t,
				\end{cases}\label{Riccati 4-4}\\
				\check{P}_i(t,s)&=e^{A'(s-t)}\check{P}_i(s,s)e^{A(s-t)},\label{Riccati 5-5}\\
				\hat{S}_i(t)&:=P_i(t)+\int_{0}^{h_1}\hat{P}_i(t,t+\theta)d\theta+\int_{0}^{h_2}\check{P}_i(t,t+\theta)d\theta,\label{Riccati 6-6}\\
				\check{S}_i(t)&:=P_i(t)+\int_{h_1-h_2}^{h_1}\hat{P}_i(t,t+\theta)d\theta+\int_{0}^{h_2}\check{P}_i(t,t+\theta)d\theta,\label{Riccati 7-7}
			\end{align}
			with $P_i(T)=H_i$, and $\hat{P}_i(T,T+\theta)=\check{P}_i(T,T+\theta)=0$ for any $\theta\geq0$, $i=1,2$.
		\end{mythm}
		\noindent\textit{Proof}. The proof consists of four steps. 

         {\it Step 1.} Firstly, with $\delta$ tending to zero, we could generate the continuous-time form of explicit solution for DFBSDEs (\ref{DFBSDEs}). Denote the limitation of $\hat{P}_{i,k}, \dfrac{1}{\delta}\hat{P}_{i,k}^{k+j}$ and $\dfrac{1}{\delta}\check{P}_{i,k}^{k+j}$ as $P_i(t), \hat{P}_i(t,t+\theta), \check{P}_i(t,t+\theta)$, respectively. 
		
		Taking the limit on the both sides of (\ref{p k explicit}), we could obtain (\ref{p i}) as
		\begin{align*}
			p_i(t)&=P_i(t)x(t)+\int_{0}^{h_1}\hat{P}_i(t,t+\theta)\mathbb{E}_{t-h_1+\theta}[x(t)]d\theta+\int_{0}^{h_2}\check{P}_i(t,t+\theta)\mathbb{E}_{t-h_2+\theta}[x(t)]d\theta.
		\end{align*} 
		Recalling (\ref{p k explicit}), we have
		\begin{align*}
			q_{i,k}&=\frac{1}{\delta}\hat{\mathbb{E}}_{k-1}(\Delta w_kp_{i,k})= \frac{1}{\delta}\hat{\mathbb{E}}_{k-1}\bigg\{\Delta w_k\hat{P}_{i,k+1}\big[(\hat{A}+\Delta w_k\bar{A})x_k
             +(\hat{B}_{11}+\Delta w_k\bar{B}_{11})\hat{\mathbb{E}}_{k-d_1-1}(p_{1,k})\\
         		&\quad +(\hat{B}_{12}+\Delta w_k\bar{B}_{12})\hat{\mathbb{E}}_{k-d_2-1}(p_{2,k})+(\hat{B}_{21}+\Delta w_k\bar{B}_{21})\hat{\mathbb{E}}_{k-d_1-1}(q_{1,k})+(\hat{B}_{22}+\Delta w_k\bar{B}_{22})\\
         		&\quad \times \hat{\mathbb{E}}_{k-d_2-1}(q_{2,k})\big]+\Delta w_k\sum_{j=0}^{d_1-1}\hat{P}_{i,k+1}^{k+1+j}\hat{\mathbb{E}}_{k-d_1+j}(x_{k+1})+\Delta w_k\sum_{j=0}^{d_2-1}\check{P}_{i,k+1}^{k+1+j}\hat{\mathbb{E}}_{k-d_2+j}(x_{k+1})\bigg\}
            \end{align*}
            \begin{align*}
				&=\hat{P}_{i,k+1}\big[\bar{A}x_k+\bar{B}_{11}\hat{\mathbb{E}}_{k-d_1-1}(p_{1,k})+\bar{B}_{12}\hat{\mathbb{E}}_{k-d_2-1}(p_{2,k}) +\bar{B}_{21}\hat{\mathbb{E}}_{k-d_1-1}(q_{1,k})\\
				&\quad+\bar{B}_{22}\hat{\mathbb{E}}_{k-d_2-1}(q_{2,k})\big].&
		\end{align*} 
		Taking the limit on the both sides, (\ref{q i}) could be obtained. Then considering the definition of $\hat{B}_{ij}$, and $\hat{S}_{i,k}, S_{i,k}^m, \check{S}_{i,k}, i,j=1,2$, we verify (\ref{x explicit}). By computation, the inverse of matrices $\hat{\Gamma}_k, \Gamma_k^m, \check{\Gamma}_k$ are as follows:
		\begin{align*}
			(\hat{\Gamma}_k)^{-1}&=\left[\begin{array}{cc}
				I-\delta B_{11}\hat{S}_{1,k}-\delta B_{12}\hat{S}_{2,k} & -(B_{21}\hat{P}_{1,k}+B_{22}\hat{P}_{2,k}) \\ -\delta(\bar{B}_{11}\hat{S}_{1,k}+\bar{B}_{12}\hat{S}_{2,k}) & I-\bar{B}_{21}\hat{P}_{1,k}-\bar{B}_{22}\hat{P}_{2,k} 
				\end{array}\right]^{-1}\equiv\left[\begin{array}{cc}
					\hat{N}_1 & \hat{N}_2 \\ \hat{N}_3 & \hat{N}_4
				\end{array}\right],\\
			(\Gamma_k^m)^{-1}&=\left[\begin{array}{cc}
				I-\delta B_{12}S_{2,k}^m & 		-B_{22}\hat{P}_{2,k} \\ -\delta\bar{B}_{12}S_{2,k}^m & I-\bar{B}_{22}\hat{P}_{2,k} 
			\end{array}\right]^{-1}\equiv\left[\begin{array}{cc}
				N_1^m & N_2^m \\ N_3^m & N_4^m
			\end{array}\right],\\
			(\check{\Gamma}_k)^{-1}&=\left[\begin{array}{cc}
				I-\delta B_{12}\check{
				S}_{2,k} & -B_{22}\hat{P}_{2,k} \\ -\delta\bar{B}_{12}\check{S}_{2,k} & I-\bar{B}_{22}\hat{P}_{2,k} 
				\end{array}\right]^{-1}\equiv\left[\begin{array}{cc}
					\check{N}_1 & \check{N}_2 \\ \check{N}_3 & \check{N}_4
				\end{array}\right],
		\end{align*}
		where
		\begin{align*}
			\hat{N}_1&:=\Big\{I-\delta(B_{11}\hat{S}_{1,k}+B_{12}\hat{S}_{2,k})-\delta(B_{21}\hat{P}_{1,k}+B_{22}\hat{P}_{2,k})\\
			         &\qquad \times(I-\bar{B}_{21}\hat{P}_{1,k}-\bar{B}_{22}\hat{P}_{2,k})^{-1}(\bar{B}_{11}\hat{S}_{1,k}+\bar{B}_{12}\hat{S}_{2,k})\Big\}^{-1},\\
			\hat{N}_2&:=\hat{N}_1(B_{21}\hat{P}_{1,k}+B_{22}\hat{P}_{2,k})(I-\bar{B}_{21}\hat{P}_{1,k}-\bar{B}_{22}\hat{P}_{2,k})^{-1},\\
			\hat{N}_3&:=\delta(I-\bar{B}_{21}\hat{P}_{1,k}-\bar{B}_{22}\hat{P}_{2,k})^{-1}(\bar{B}_{11}\hat{S}_{1,k}+\bar{B}_{12}\hat{S}_{2,k})\hat{N}_1,\\
			\hat{N}_4&:=\big[I+\hat{N}_3(B_{21}\hat{P}_{1,k}+B_{22}\hat{P}_{2,k})\big](I-\bar{B}_{21}\hat{P}_{1,k}-\bar{B}_{22}\hat{P}_{2,k})^{-1},\\
			N_1^m&:=\big[I-\delta B_{12}S_{2,k}^m-\delta B_{22}\hat{P}_{2,k}(I-\bar{B}_{22}\hat{P}_{2,k})^{-1}\bar{B}_{12}S_{2,k}^m\big]^{-1},\\
			N_2^m&:=N_1^mB_{22}\hat{P}_{2,k}(I-\bar{B}_{22}\hat{P}_{2,k})^{-1},\quad N_3^m:=\delta(I-\bar{B}_{22}\hat{P}_{2,k})^{-1}\bar{B}_{12}S_{2,k}^mN_1^m,\\
			N_4^m&:=(I+N_3^mB_{22}\hat{P}_{2,k})(I-\bar{B}_{22}\hat{P}_{2,k})^{-1},\\
			\check{N}_1&:=\big[I-\delta B_{12}\check{S}_{2,k}-\delta B_{22}\hat{P}_{2,k}(I-\bar{B}_{22}\hat{P}_{2,k})^{-1}\bar{B}_{12}\check{S}_{2,k}\big]^{-1},\\
			\check{N}_2&:=\check{N}_1B_{22}\hat{P}_{2,k}(I-\bar{B}_{22}\hat{P}_{2,k})^{-1},\quad \check{N}_3:=\delta(I-\bar{B}_{22}\hat{P}_{2,k})^{-1}\bar{B}_{12}\check{S}_{2,k}\check{N}_1,\\
			\check{N}_4&:=(I+\check{N}_3B_{22}\hat{P}_{2,k})(I-\bar{B}_{22}\hat{P}_{2,k})^{-1}.
		\end{align*}
		Noting the form of coefficient matrices, (\ref{x k explicit}) can be rewritten as follows,
		\begin{align*}
			x_k&=x_{k-1}+\delta\Big\{Ax_{k-1}+\big\{\big[B_{11}+B_{22}\hat{P}_{2,k}(I-\bar{B}_{22}\hat{P}_{2,k})^{-1}\bar{B}_{11}\big]\hat{S}_{1,k}\hat{N}_1+\big[B_{21}+B_{22}\hat{P}_{2,k}\\
			&\qquad \times(I-\bar{B}_{22}\hat{P}_{2,k})^{-1}\bar{B}_{21}\big]\hat{P}_{1,k}(I-\bar{B}_{21}\hat{P}_{1,k}-\bar{B}_{22}\hat{P}_{2,k})^{-1}\big[(\bar{B}_{11}\hat{S}_{1,k}+\bar{B}_{12}\hat{S}_{2,k})\hat{N}_1+\bar{A}\big]\big\}\\
			&\qquad \times\hat{\mathbb{E}}_{k-d_1-2}(x_{k-1})+\sum_{j=1}^{d_1-d_2-1}\big[B_{12}+B_{22}\hat{P}_{2,k}(I-\bar{B}_{22}\hat{P}_{2,k})^{-1}\bar{B}_{12}\big]\hat{P}_{2,k}^{k-1+j}N_1^j\hat{\mathbb{E}}_{k-d_1-2+j}(x_{k-1})\\
			&\qquad +\big[B_{12}\check{S}_{2,k}\check{N}_1+B_{22}\hat{P}_{2,k}(I-\bar{B}_{22}\hat{P}_{2,k})^{-1}(\bar{B}_{12}\check{S}_{2,k}\check{N}_1+\bar{A})\big]\hat{\mathbb{E}}_{k-d_2-2}(x_{k-1})\Big\}
		\end{align*}
		\begin{align*}
			&\quad +\Delta w_k\Big\{\bar{A}x_{k-1}+\big\{\big[\bar{B}_{11}+\bar{B}_{22}\hat{P}_{2,k}(I-\bar{B}_{22}\hat{P}_{2,k})^{-1}\bar{B}_{11}\big]\hat{S}_{1,k}\hat{N}_1+\big[\bar{B}_{21}+\bar{B}_{22}\hat{P}_{2,k}\\
			&\qquad \times(I-\bar{B}_{22}\hat{P}_{2,k})^{-1}\bar{B}_{21}\big]\hat{P}_{1,k}(I-\bar{B}_{21}\hat{P}_{1,k}-\bar{B}_{22}\hat{P}_{2,k})^{-1}\big[(\bar{B}_{11}\hat{S}_{1,k}+\bar{B}_{12}\hat{S}_{2,k})\hat{N}_1+\bar{A}\big]\big\}\\
			&\qquad \times\hat{\mathbb{E}}_{k-d_1-2}(x_{k-1})+\sum_{j=1}^{d_1-d_2-1}\big[\bar{B}_{12}+\bar{B}_{22}\hat{P}_{2,k}(I-\bar{B}_{22}\hat{P}_{2,K})^{-1}\bar{B}_{12}\big]\hat{P}_{2,k}^{k-1+j}N_1^j\hat{\mathbb{E}}_{k-d_1-2+j}(x_{k-1})\\
			&\qquad +\big[\bar{B}_{12}\check{S}_{2,k}\check{N}_1+\bar{B}_{22}\hat{P}_{2,k}(I-\bar{B}_{22}\hat{P}_{2,k})^{-1}(\bar{B}_{12}\check{S}_{2,k}\check{N}_1+\bar{A})\big]\hat{\mathbb{E}}_{k-d_2-2}(x_{k-1})\Big\}+o(\delta).
		\end{align*}
		Letting $\delta\equiv t_{k+1}-t_k \rightarrow 0$, we have (\ref{x explicit}). 
		
		{\it Step 2.} Secondly, we verify (\ref{Riccati 1-1})-(\ref{Riccati 7-7}) which are the system of Riccati equations. Inputting $\hat{A}':=I+\delta A', \hat{A}:=I+\delta A$ and $\hat{Q}_i:=\delta Q_i$ into (\ref{Riccati 1}), we get
		\begin{align*}
			\hat{P}_{i,k}&=\hat{A}'\hat{P}_{i,k+1}\hat{A}+\delta\bar{A}'\hat{P}_{i,k+1}\bar{A}+\hat{A}'(\hat{P}_{i,k+1}^{k+1+d_1}+\check{P}_{i,k+1}^{k+1+d_2})\hat{A}+\hat{Q}_i\\
			&=(I+\delta A')\hat{P}_{i,k+1}(I+\delta A)+\delta\bar{A}'\hat{P}_{i,k+1}\bar{A}+(I+\delta A')(\hat{P}_{i,k+1}^{k+1+d_1}+\check{P}_{i,k+1}^{k+1+d_2})(I+\delta A)+\delta Q_i.
		\end{align*}
		Dividing both sides by $\delta$, it yields
		\begin{align*}
			\frac{\hat{P}_{i,k}-\hat{P}_{i,k+1}}{\delta}=A'\hat{P}_{i,k+1}+\hat{P}_{i,k+1}A+\bar{A}'\hat{P}_{i,k+1}\bar{A}+\frac{1}{\delta}(I+\delta A')(\hat{P}_{i,k+1}^{k+1+d_1}+\check{P}_{i,k+1}^{k+1+d_2})(I+\delta A)+Q_i.
		\end{align*}
		Letting $\delta \rightarrow 0$, we could get
		\begin{align*}
			-\dot{P}_i(t)=A'P_i(t)+P_i(t)A+\bar{A}'\hat{P}_i(t)\bar{A}+\hat{P}_i(t,t+h_1)+\check{P}_i(t,t+h_2)+Q_i.
		\end{align*}
		Thus we obtain (\ref{Riccati 1-1}). Then we prove  (\ref{Riccati 3-3}) and  (\ref{Riccati 4-4}). We have
		\begin{align*}
			\check{P}_{i,k}^k=\mathbb{E}\big[(\hat{A}'\check{S}_{i,k+1}+\Delta w_k\bar{A}'\hat{P}_{i,k+1})H_k\big],
		\end{align*}
		where 
		\begin{align*}
			H_k&:=\delta\big[B_{12}\check S_{2,k+1}\check{N}_1+B_{22}\hat{P}_{2,k+1}(I-\bar{B}_{22}\hat{P}_{2,k+1})^{-1}(\bar{B}_{12}\check S_{2,k+1}\check{N}_1+\bar{A})\big]\\
			&\quad +\Delta w_k\big[\bar{B}_{12}\check S_{2,k+1}\check{N}_1+\bar{B}_{22}\hat{P}_{2,k+1}(I-\bar{B}_{22}\hat{P}_{2,k+1})^{-1}(\bar{B}_{12}\check S_{2,k+1}\check{N}_1+\bar{A})\big]+o(\delta).
		\end{align*}
		Dividing both sides by $\delta$ and letting $\delta \rightarrow 0$, we achieve
		\begin{align*}
			\check{P}_i(t,t)&=\check{S}_i(t)\big[B_{12}\check{S}_2(t)+B_{22}P_2(t)(I-\bar{B}_{22}P_2(t))^{-1}(\bar{B}_{12}\check{S}_2(t)+\bar{A})\big]\\
            &\quad +\bar{A}'P_i(t)\big[\bar{B}_{12}\check{S}_2(t)+\bar{B}_{22}P_2(t)(I-\bar{B}_{22}P_2(t))^{-1}(\bar{B}_{12}\check{S}_2(t)+\bar{A})\big]\\
			&=[\check{S}_i(t)B_{12}+\bar{A}'P_i(t)\bar{B}_{12}]\check{S}_2(t)+[\check{S}_i(t)B_{22}+\bar{A}'P_i(t)\bar{B}_{22}]P_2(t)\\
			&\quad \times(I-\bar{B}_{22}P_2(t))^{-1}(\bar{B}_{12}\check{S}_2(t)+\bar{A}),
		\end{align*}
		 we prove (\ref{Riccati 3-3}). Recalling the definition of $\hat{A}'$ and $\hat{A}$, $\check{P}_{i,k}^{k+j}$ can be rewritten as 
		 \begin{align*}
		 	\check{P}_{i,k}^{k+j}=\hat{A}'\check{P}_{i,k+1}^{k+j}\hat{A}=\check{P}_{i,k+1}^{k+j}+\delta A'\check{P}_{i,k+1}^{k+j}+\delta\check{P}_{i,k+1}^{k+j}A+\delta^2A'\check{P}_{i,k+1}^{k+j}A.
		 \end{align*}
		Similarly to the deduction of $P_i(t)$, it yields
		\begin{align*}
			-\partial_t\check{P}(t,s)=A'\check{P}(t,s)+\check{P}(t,s)A.
		\end{align*}
		Solving the equation, it generates $\check{P}_i(t,s)=e^{A'(s-t)}\check{P}(s,s)e^{A(s-t)}$. The deductions of $\hat{P}(t,t)$ and $\hat{P}(t,s)$ have no difference, so we omitted it here. So far, we obtained (\ref{Riccati 1-1})-(\ref{Riccati 5-5}).

        {\it Step 3.} Next, by using It\^{o}'s formula, we can verify that the continuous-time form of explicit solution is exactly the solution of DFBSDEs (\ref{DFBSDEs}). We denote
		\begin{align}
			\hat{p}_i(t)&:=P_i(t)x(t)+\int_{0}^{h_1}\hat{P}_i(t,t+\theta)\mathbb{E}_{t-h_1+\theta}[x(t)]d\theta+\int_{0}^{h_2}\check{P}_i(t,t+\theta)\mathbb{E}_{t-h_2+\theta}[x(t)]d\theta,\label{hat p i}\\
			\hat{q}_i(t)&:=P_i(t)\Big\{\bar{A}x(t)+\bar{B}_{11}\mathbb{E}_{t-h_1}[\hat{p}_1(t)]+\bar{B}_{12}\mathbb{E}_{t-h_2}[\hat{p}_2(t)]\nonumber\\
            &\qquad +\bar{B}_{21}\mathbb{E}_{t-h_1}[\hat{q}_1(t)]+\bar{B}_{22}\mathbb{E}_{t-h_2}[\hat{q}_2(t)]\Big\}\label{hat q i}.
		\end{align}
		Taking conditional expectation on both sides of (\ref{hat p i}) and (\ref{hat q i}), they generate
		\begin{align*}
			\mathbb{E}_{t-h_1}[\hat{p}_i(t)]&=\bigg[P_i(t)+\int_0^{h_1}\hat{P}_i(t,t+\theta)d\theta+\int_{0}^{h_2}\check{P}_i(t,t+\theta)d\theta\bigg]\mathbb{E}_{t-h_1}[x(t)]\\
			&:=\hat{S}_i(t)\mathbb{E}_{t-h_1}[x(t)],\\
			\mathbb{E}_{t-h_2}[\hat{p}_2(t)]&=\bigg[P_2(t)+\int_{h_1-h_2}^{h_1}\hat{P}_2(t,t+\theta)d\theta+\int_{0}^{h_2}\check{P}_2(t,t+\theta)d\theta\bigg]\mathbb{E}_{t-h_2}[x(t)]\\
			&\quad+\int_0^{h_1-h_2}\hat{P}_2(t,t+\theta)\mathbb{E}_{t-h_1+\theta}[x(t)]d\theta\\
			&:=\check{S}_2(t)\mathbb{E}_{t-h_2}[x(t)]+\int_{0}^{h_1-h_2}\hat{P}_2(t,t+\theta)\mathbb{E}_{t-h_1+\theta}[x(t)]d\theta,\\
			\mathbb{E}_{t-h_1}[\hat{q}_1(t)]&=(I-\bar{B}_{21}P_1(t)-\bar{B}_{22}P_2(t))^{-1}P_1(t)\big\{\bar{A}\mathbb{E}_{t-h_1}[x(t)]+\bar{B}_{11}\mathbb{E}_{t-h_1}[\hat{p}_1(t)]\\
			&\quad+\bar{B}_{12}\mathbb{E}_{t-h_1}[\hat{p}_2(t)]\big\},\\	
            \mathbb{E}_{t-h_2}[\hat{q}_2(t)]&=(I-\bar{B}_{22}P_2(t))^{-1}P_2(t)\big\{\bar{B}_{21}(I-\bar{B}_{21}P_1(t)-\bar{B}_{22}P_2(t))^{-1}P_1(t)\bar{A}\mathbb{E}_{t-h_1}[x(t)]\\
			&\quad+\bar{A}\mathbb{E}_{t-h_2}[x(t)]+[\bar{B}_{11}+\bar{B}_{21}(I-\bar{B}_{21}P_1(t)-\bar{B}_{22}P_2(t))^{-1}P_1(t)\bar{B}_{11}]\mathbb{E}_{t-h_1}[\hat{p}_1(t)]\\
			&\quad+\bar{B}_{21}(I-\bar{B}_{21}P_1(t)-\bar{B}_{22}P_2(t))^{-1}P_1(t)\bar{B}_{12}\mathbb{E}_{t-h_1}[\hat{p}_2(t)]+\bar{B}_{12}\mathbb{E}_{t-h_2}[\hat{p}_2(t)]\big\}.
			\end{align*}
		Considering the definition of $\hat{S}_i(t), \check{S}_i(t)$, we rewrite (\ref{x explicit}). The $dt$ term of the right hand of (\ref{x explicit}) can be rewritten as follows,
		\begin{align*}
			&Ax(t)+B_{11}\mathbb{E}_{t-h_1}[\hat{p}_1(t)]+B_{12}\mathbb{E}_{t-h_2}[\hat{p}_2(t)]+B_{21}P_1(t)(I-\bar{B}_{21}P_1(t)-\bar{B}_{22}P_2(t))^{-1}\\
			&\quad \times\big\{\bar{A}\mathbb{E}_{t-h_1}[x(t)]+\bar{B}_{11}\mathbb{E}_{t-h_1}[\hat{p}_1(t)]+\bar{B}_{12}\mathbb{E}_{t-h_2}[\hat{p}_2(t)]\big\}+B_{22}P_2(t)(I-\bar{B}_{22}P_2(t))^{-1}\\
			&\quad \times\big\{\bar{B}_{21}(I-\bar{B}_{21}P_1(t)-\bar{B}_{22}P_2(t))^{-1}P_1(t)\bar{A}\mathbb{E}_{t-h_1}[x(t)]+\bar{A}\mathbb{E}_{t-h_2}[x(t)]+[\bar{B}_{11}+\bar{B}_{21}\\
			&\quad \times(I-\bar{B}_{21}P_1(t)-\bar{B}_{22}P_2(t))^{-1}P_1(t)\bar{B}_{11}]\mathbb{E}_{t-h_1}[\hat{p}_1(t)]+\bar{B}_{21}(I-\bar{B}_{21}P_1(t)-\bar{B}_{22}P_2(t))^{-1}\\
			&\quad \times P_1(t)\bar{B}_{12}\mathbb{E}_{t-h_1}[\hat{p}_2(t)]+\bar{B}_{12}\mathbb{E}_{t-h_2}[\hat{p}_2(t)]\big\}\\
			&=Ax(t)+B_{11}\mathbb{E}_{t-h_1}[\hat{p}_1(t)]+B_{12}\mathbb{E}_{t-h_2}[\hat{p}_2(t)]+B_{21}\mathbb{E}_{t-h_1}[\hat{q}_1(t)]+B_{22}\mathbb{E}_{t-h_2}[\hat{q}_2(t)].		
		\end{align*}
		With the same procedure, we can also deal with the term of $dw(t)$. Then $x(t)$ satisfies the equality as follows,
		\begin{align}
			dx(t)=&\big\{Ax(t)+B_{11}\mathbb{E}_{t-h_1}[\hat{p}_1(t)]+B_{12}\mathbb{E}_{t-h_2}[\hat{p}_2(t)]+B_{21}\mathbb{E}_{t-h_1}[\hat{q}_1(t)]\nonumber\\
			&+B_{22}\mathbb{E}_{t-h_2}[\hat{q}_2(t)]\big\}dt+\big\{\bar{A}x(t)+\bar{B}_{11}\mathbb{E}_{t-h_1}[\hat{p}_1(t)]+\bar{B}_{12}\mathbb{E}_{t-h_2}[\hat{p}_2(t)]\nonumber\\
			&+\bar{B}_{21}\mathbb{E}_{t-h_1}[\hat{q}_1(t)]+\bar{B}_{22}\mathbb{E}_{t-h_2}[\hat{q}_2(t)]\big\}dw(t).\label{d x(t)}
		\end{align}
		Considering (\ref{d x(t)}), (\ref{Riccati 1-1}) - (\ref{Riccati 5-5}), utilizing It\^{o}'s formula, we can get the following results for $i=1$ while $i=2$ can be derived in the same way. In fact, we have
 		\begin{align}
			d[P_1(t)x(t)]=&-\big[A'P_1(t)+P_1(t)A+\bar{A}'P_1(t)\bar{A}+Q_1+\hat{P}_1(t,t+h_1)+\check{P}_1(t,t+h_2)\big]x(t)dt\nonumber\\
				&+P_1(t)\big\{Ax(t)+B_{11}\mathbb{E}_{t-h_1}[\hat{p}_1(t)]+B_{12}\mathbb{E}_{t-h_2}[\hat{p}_2(t)]+B_{21}\mathbb{E}_{t-h_1}[\hat{q}_1(t)]\nonumber\\
				&+B_{22}\mathbb{E}_{t-h_2}[\hat{q}_2(t)]\big\}dt+P_1(t)\big\{\bar{A}x(t)+\bar{B}_{11}\mathbb{E}_{t-h_1}[\hat{p}_1(t)]+\bar{B}_{12}\mathbb{E}_{t-h_2}[\hat{p}_2(t)]\nonumber\\
				&+\bar{B}_{21}\mathbb{E}_{t-h_1}[\hat{q}_1(t)]+\bar{B}_{22}\mathbb{E}_{t-h_2}[\hat{q}_2(t)]\big\}dw(t).\label{Ito formula}
		\end{align}
		Noting that
		\begin{align*}
			&d\bigg[\int_0^{h_1}\hat{P}_1(t,t+\theta)\mathbb{E}_{t-h_1+\theta}[x(t)]d\theta\bigg]=\int_t^{t+h_1}\partial_t\hat{P}_1(t,s)\mathbb{E}_{s-h_1}[x(t)]ds\\
			&+\int_{t}^{t+h_1}\hat{P}_1(t,s)\partial_t\mathbb{E}_{s-h_1}[x(t)]ds+\hat{P}_1(t,t+h_1)x(t)dt-\hat{P}_1(t,t)\mathbb{E}_{t-h_1}[x(t)]dt,
		\end{align*}
		where $\partial_t\hat{P}_1(t,s)$ satisfying (\ref{Riccati 4-4}). Thus it generates
		\begin{align}
			&d\bigg[\int_0^{h_1}\hat{P}_1(t,t+\theta)\mathbb{E}_{t-h_1+\theta}[x(t)]d\theta\bigg]=-A'\int_{t+h_1-h_2}^{t+h_1}\partial_t\hat{P}_1(t,s)\mathbb{E}_{s-h_1}[x(t)]dsdt\nonumber\\
			&+\int_t^{t+h_1-h_2}\partial_t\hat{P}_1(t,s)\mathbb{E}_{s-h_1}[x(t)]ds+\int_t^{t+h_1}\hat{P}_1(t,s)ds\big\{B_{11}\mathbb{E}_{t-h_1}[\hat{p}_1(t)]\nonumber\\
			&+B_{21}\mathbb{E}_{t-h_1}[\hat{q}_1(t)]\big\}dt+\int_t^{t+h_1-h_2}\hat{P}_1(t,s)\big\{B_{12}\mathbb{E}_{s-h_1}[\hat{p}_2(t)]+B_{22}\mathbb{E}_{s-h_1}[\hat{q}_2(t)]\big\}dsdt\nonumber\\
			&+\int_{t+h_1-h_2}^{t+h_1}\hat{P}_1(t,s)ds\big\{B_{12}\mathbb{E}_{t-h_2}[\hat{p}_2(t)]+B_{22}\mathbb{E}_{t-h_2}[\hat{q}_2(t)]\big\}dt\nonumber\\
			&+\hat{P}_1(t,t+h_1)x(t)dt-\hat{P}_1(t,t)\mathbb{E}_{t-h_1}[x(t)]dt.\label{d P 1 hat}
		\end{align}
		In the same way, we can also derive
		\begin{align}
			&d\bigg[\int_0^{h_2}\check{P}_1(t,t+\theta)\mathbb{E}_{t-h_2+\theta}[x(t)]d\theta\bigg]=-A'\int_t^{t+h_2}\check{P}_1(t,s)\mathbb{E}_{s-h_2}[x(t)]ds\nonumber\\
			&+\int_t^{t+h_2}\check{P}_1(t,s)ds\big\{B_{11}\mathbb{E}_{t-h_1}[\hat{p}_1(t)]+B_{12}\mathbb{E}_{t-h_2}[\hat{p}_2(t)]+B_{21}\mathbb{E}_{t-h_1}[\hat{q}_1(t)]\nonumber\\
			&+B_{22}\mathbb{E}_{t-h_2}[\hat{q}_2(t)]\big\}+\check{P}_1(t,t+h_2)x(t)dt-\check{P}_1(t,t)\mathbb{E}_{t-h_2}[x(t)]dt.\label{d P 1 check}
		\end{align}
		Combining (\ref{Ito formula}), (\ref{d P 1 hat}) and (\ref{d P 1 check}), the following result can be derived:
		\begin{align*}
			d\bigg[&P_1(t)x(t)+\int_{0}^{h_1}\hat{P}_1(t,t+\theta)\mathbb{E}_{t-h_1+\theta}[x(t)]d\theta+\int_{0}^{h_2}\check{P}_1(t,t+\theta)\mathbb{E}_{t-h_2+\theta}[x(t)]d\theta\bigg]\\
			=&-\bigg\{A'[P_1(t)x(t)+\int_{0}^{h_1}\hat{P}_1(t,t+\theta)\mathbb{E}_{t-h_1+\theta}[x(t)]d\theta+\int_{0}^{h_2}\check{P}_1(t,t+\theta)\mathbb{E}_{t-h_2+\theta}[x(t)]d\theta]\bigg\}dt\\
			&-A'P_1(t)\Big\{\bar{A}x(t)+\bar{B}_{11}\mathbb{E}_{t-h_1}[\hat{p}_1(t)]+\bar{B}_{12}\mathbb{E}_{t-h_2}[\hat{p}_2(t)]+\bar{B}_{21}\mathbb{E}_{t-h_1}[\hat{q}_1(t)]\\
			&+\bar{B}_{22}\mathbb{E}_{t-h_2}[\hat{q}_2(t)]\Big\}dt-Q_1x(t)dt+P_1(t)\Big\{\bar{A}x(t)+\bar{B}_{11}\mathbb{E}_{t-h_1}[\hat{p}_1(t)]+\bar{B}_{12}\mathbb{E}_{t-h_2}[\hat{p}_2(t)]\\
			&+\bar{B}_{21}\mathbb{E}_{t-h_1}[\hat{q}_1(t)]+\bar{B}_{22}\mathbb{E}_{t-h_2}[\hat{q}_2(t)]\Big\}dw(t)+\Big\{[P_1(t)B_{11}+\bar{A}'P_1(t)\bar{B}_{11}]\mathbb{E}_{t-h_1}[\hat{p}_1(t)]\\
			&+[P_1(t)B_{12}+\bar{A}'P_1(t)\bar{B}_{12}]\mathbb{E}_{t-h_2}[\hat{p}_2(t)]+[P_1(t)B_{21}+\bar{A}'P_1(t)\bar{B}_{21}]\mathbb{E}_{t-h_1}[\hat{q}_1(t)]\\
            &+[P_1(t)B_{22}+\bar{A}'P_1(t)\bar{B}_{22}]\mathbb{E}_{t-h_2}[\hat{q}_2(t)]\Big\}dt+\int_t^{t+h_1-h_2}\partial_t\hat{P}_1(t,s)\mathbb{E}_{s-h_1}[x(t)]ds\\
			&+\int_t^{t+h_1}\hat{P}_1(t,s)ds\big\{B_{11}\mathbb{E}_{t-h_1}[\hat{p}_1(t)]+B_{21}\mathbb{E}_{t-h_1}[\hat{q}_1(t)]\big\}dt\\
			&+\int_t^{t+h_1-h_2}\hat{P}_1(t,s)\big\{B_{12}\mathbb{E}_{s-h_1}[\hat{p}_2(t)]+B_{22}\mathbb{E}_{s-h_1}[\hat{q}_2(t)]\big\}dsdt\\
			&+\int_{t+h_1-h_2}^{t+h_1}\hat{P}_1(t,s)ds\big\{B_{12}\mathbb{E}_{t-h_2}[\hat{p}_2(t)]+B_{22}\mathbb{E}_{t-h_2}[\hat{q}_2(t)]\big\}dt-\hat{P}_1(t,t)\mathbb{E}_{t-h_1}[x(t)]dt\\
			&+\int_t^{t+h_2}\check{P}_1(t,s)ds\big\{B_{11}\mathbb{E}_{t-h_1}[\hat{p}_1(t)]+B_{12}\mathbb{E}_{t-h_2}[\hat{p}_2(t)]\\
			&+B_{21}\mathbb{E}_{t-h_1}[\hat{q}_1(t)]+B_{22}\mathbb{E}_{t-h_2}[\hat{q}_2(t)]\big\}-\check{P}_1(t,t)\mathbb{E}_{t-h_2}[x(t)]dt.
		\end{align*}
		Substituting  $\hat{S}_1(t), \check{S}_1(t), \hat{P}_1(t,t)$ and $\check{P}_1(t,t)$ into the above equation, it yields 
		\begin{align}
			d\hat{p}_1(t)=-[A'\hat{p}_1(t)+\bar{A}'\hat{q}_1(t)+Q_1x(t)]dt+\hat{q}_1(t)dw(t),\label{p 1}
		\end{align}
		and 
		\begin{align}
			\hat{q}_1(t)=&P_1(t)\{\bar{A}x(t)+\bar{B}_{11}\mathbb{E}_{t-h_1}[\hat{p}_1(t)]+\bar{B}_{12}\mathbb{E}_{t-h_2}[\hat{p}_2(t)]+\bar{B}_{21}\mathbb{E}_{t-h_1}[\hat{q}_1(t)]\nonumber\\
			&+\bar{B}_{22}\mathbb{E}_{t-h_2}[\hat{q}_2(t)]\}.\label{q 1}
		\end{align}
		By comparison, it can be seen that (\ref{d x(t)}), (\ref{p 1}) and (\ref{q 1}) are same as the forward and backward equation of DFBSDEs (\ref{DFBSDEs}). Therefore, (\ref{p i})-(\ref{x explicit}) is the solution of (\ref{DFBSDEs}).
				 
		{\it Step 4.} Finally, we give the proof for the uniqueness of the explicit solution. Let $\zeta_i(t):=P_i(t)x(t)+\int_0^{h_1}\hat{P}_i(t,t+\theta)\mathbb{E}_{t-h_1+\theta}[x(t)]d\theta+\int_0^{h_2}\check{P}_i(t,t+\theta)\mathbb{E}_{t-h_2+\theta}[x(t)]d\theta+\eta_i(t)$, for some $\eta_i(\cdot)$ with $\eta_i(T)=0, i=1,2$. Using It\^{o}'s formula, we have
		\begin{align*}
			d\eta_1(t)&=d\zeta_1(t)-d\bigg[P_1(t)x(t)+\int_0^{h_1}\hat{P}_1(t,t+\theta)\mathbb{E}_{t-h_1+\theta}[x(t)]d\theta\\
            &\quad +\int_0^{h_2}\check{P}_1(t,t+\theta)\mathbb{E}_{t-h_2+\theta}[x(t)]d\theta\bigg]
        \end{align*}
        \begin{align*}
            &=-\big[A'\zeta_1(t)+\bar{A}'r_1(t)+Q_1x(t)\big]dt+r_1(t)dw(t)\\
            &\quad +\bigg\{A'\bigg[P_1(t)x(t)+\int_{0}^{h_1}\hat{P}_1(t,t+\theta)\mathbb{E}_{t-h_1+\theta}[x(t)]d\theta+\int_0^{h_2}\check{P}_1(t,t+\theta)\mathbb{E}_{t-h_2+\theta}[x(t)]d\theta\bigg]\\
            &\quad +\bar{A}'q_1(t)+Q_1x(t)\bigg\}dt-q_1(t)dw(t)\\
			&=-\big\{A'\eta_1(t)+\bar{A}'[r_1(t)-q_1(t)]\big\}dt+[r_1(t)-q_1(t)]dw(t),
		\end{align*}     
		for $i=1$. With $\eta_1(T)=0$, it implies $\eta_1(t)=0$, while $i=2$ has no difference. Thus (\ref{p i})-(\ref{x explicit}) is the unique solution.$\hfill\qedsymbol$

		\begin{Remark}
			It must be pointed out that	the derivation of the system of Riccati equations (\ref{Riccati 1-1})-(\ref{Riccati 5-5}) is formal. However, we have provided a rigorous proof which verified that the continuous-time form of explicit solution indeed satisfy DFBSDEs (\ref{DFBSDEs}).
		\end{Remark}

		In the end, the state-estimate feedback Nash equilibrium of \textbf{Problem SDG-ADI} is presented in the below theorem.
		\begin{mythm}
			If Riccati equations (\ref{Riccati 1-1})-(\ref{Riccati 5-5}) admit solutions such that matrices $R_1(\cdot), R_2(\cdot)$ are invertible, the Nash equilibrium of \textbf{Problem SDG-ADI} has the state-estimate feedback form as follows:
			\begin{align}
				u_1(t)&=-R_1^{-1}(t)O_1(t)\mathbb{E}_{t-h_1}[x(t)],\label{u 1}\\
				u_2(t)&=-R_2^{-1}(t)\bigg\{-\bar{B}'_2P_2(t)\bar{B}_1R_1^{-1}(t)O_1(t)\mathbb{E}_{t-h_1}[x(t)]+B'_2\int_0^{h_1-h_2}\hat{P}_2(t,t+\theta)\nonumber\\
				&\qquad \times\mathbb{E}_{t-h_1+\theta}[x(t)]d\theta+\big[B'_2\check{S}_2(t)+\bar{B}'_2P_2(t)\bar{A}\big]\mathbb{E}_{t-h_2}[x(t)]\bigg\},\label{u 2}
			\end{align}
			where
			\begin{align}
				R_1(t)&:=R_1+\bar{B}'_1P_1(t)\bar{B}_1-[\bar{B}'_1P_1(t)\bar{B}_2]R_2^{-1}(t)[\bar{B}'_2P_2(t)\bar{B}_1],\label{R 1}\\
				R_2(t)&:=R_2+\bar{B}'_2P_2(t)\bar{B}_2,\label{R 2}\\
				O_1(t)&:=B'_1\hat{S}_1(t)+\bar{B}'_1P_1(t)\bar{A}-[\bar{B}'_1P_1(t)\bar{B}_2]R_2^{-1}(t)\bigg[B'_2\check{S}_2(t)\nonumber\\
				&\quad +\bar{B}'_2P_2(t)\bar{A}+B'_2\int_0^{h_1-h_2}\hat{P}_2(t,t+\theta)d\theta\bigg].\label{O 1}
			\end{align}
		\end{mythm}

		\noindent\textit{Proof}. With $u_i(t)=-R_i^{-1}\mathbb{E}[B'_ip_i(t)+\bar{B}'_iq_i(t)|\mathcal{F}_{t-h_i}], i=1,2$ and the definition of $\bar{B}_{ij}, i,j=1,2$, we can derive
		\begin{align}
			q_i(t)&=P_i(t)\big\{\bar{A}x(t)-\bar{B}_1R_1^{-1}B'_1\mathbb{E}_{t-h_1}[\hat{p}_1(t)]-\bar{B}_2R_2^{-1}B'_2\mathbb{E}_{t-h_2}[\hat{p}_2(t)]\nonumber\\
			&\quad -\bar{B}_1R_1^{-1}\bar{B}'_1\mathbb{E}_{t-h_1}[\hat{q}_1(t)]-\bar{B}_2R_2^{-1}\bar{B}'_2\mathbb{E}_{t-h_2}[\hat{q}_2(t)]\big\}\nonumber\\
			&=P_i(t)\big[\bar{A}x(t)+\bar{B}_1u_1(t)+\bar{B}_2u_2(t)\big].\label{q i P i}
		\end{align}
		Then we deal with $u_i(t)$ by plugging (\ref{p i}), (\ref{q i P i}) into (\ref{Nash equilibrium1}) and (\ref{Nash equilibrium2}),
		\begin{align}
			0&=R_1u_1(t)+\mathbb{E}_{t-h_1}\bigg\{B'_1\bigg[P_1(t)x(t)+\int_0^{h_1}\hat{P}_1(t,t+\theta)\mathbb{E}_{t-h_1+\theta}[x(t)]d\theta\nonumber\\
			&\quad +\int_0^{h_2}\check{P}_1(t,t+\theta)\mathbb{E}_{t-h_2+\theta}[x(t)]d\theta\bigg]+\bar{B}'_1P_1(t)\big[\bar{A}x(t)+\bar{B}_1u_1(t)+\bar{B}_2u_2(t)\big]\bigg\}\nonumber\\
			&=\big[R_1+\bar{B}'_1P_1(t)\bar{B}_1\big]u_1(t)+\big[B'_1\hat{S}_1(t)+\bar{B}'_1P_1(t)\bar{A}\big]\mathbb{E}_{t-h_1}[x(t)]\nonumber\\
            &\quad +\bar{B}'_1P_1(t)\bar{B}_2\mathbb{E}_{t-h_1}[u_2(t)],\label{u 1 explicit}\\
			0&=R_2u_2(t)+\mathbb{E}_{t-h_2}\bigg\{B'_2\bigg[P_2(t)x(t)+\int_0^{h_1}\hat{P}_2(t,t+\theta)\mathbb{E}_{t-h_1+\theta}[x(t)]d\theta\nonumber\\
			&\quad +\int_{0}^{h_2}\check{P}_2(t,t+\theta)\mathbb{E}_{t-h_2+\theta}[x(t)]d\theta\bigg]+\bar{B}'_2P_2(t)\big[\bar{A}x(t)+\bar{B}_1u_1(t)+\bar{B}_2u_2(t)\big]\bigg\}\nonumber\\
			&=R_2u_2(t)+\big[B'_2\check{S}_2(t)+\bar{B}_2P_2(t)\bar{A}\big]\mathbb{E}_{t-h_2}[x(t)]+\bar{B}'_2P_2(t)\bar{B}_1u_1(t)\nonumber\\
			&\quad +\bar{B}'_2P_2(t)\bar{B}_2u_2(t)+B'_2\int_0^{h_1-h_2}\hat{P}_2(t,t+\theta)\mathbb{E}_{t-h_1+\theta}[x(t)]d\theta.\label{u 2 explicit}
		\end{align}
		(\ref{u 2 explicit}) can be rewritten as
		\begin{align}
			u_2(t)=&-\big[R_2+\bar{B}'_2P_2(t)\bar{B}_2\big]^{-1}\bigg\{\big[B'_2\check{S}_2(t)+\bar{B}_2P_2(t)\bar{A}\big]\mathbb{E}_{t-h_2}[x(t)]\nonumber\\
			&+B'_2\int_{0}^{h_1-h_2}\hat{P}_2(t,t+\theta)\mathbb{E}_{t-h_1+\theta}[x(t)]d\theta+\bar{B}'_2P_2(t)\bar{B}_1u_1(t)\bigg\}.\label{u 2 explicit -111}
		\end{align}
		By putting the above equation into (\ref{u 1 explicit}), (\ref{u 1}) can be obtained directly. Similarly, by plugging (\ref{u 1}) into (\ref{u 2 explicit -111}), we will arrive at (\ref{u 2}).$\hfill\qedsymbol$
		
\section{Concluding remarks}

	In this paper, we have studied the LQ non-zero sum differential game problem with time delays. The delays of two players are different, resulting in an asymmetric information structure of the dynamic system. To find the Nash equilibrium, we used the stochastic maximum principle to convert the solvability of the original problem into that of DFBSDEs. By using discretisation approach and backward iteration technique, we decoupled these DFBSDEs and obtained their explicit solutions. Finally, we presented the state-estimate feedback Nash equilibrium.

    It is desirable to study Stackelberg differential game with asymmetric delayed information. The partially observed system with multiple time delays, also an asymmetric information structure, is interesting but challenging. These topics will be considered in our future work.


\begin{thebibliography}{0}
		
	\bibitem{AO2014}N. Agram, B. \O ksendal, Infinite horizon optimal control of forward–backward stochastic differential equations with delay. \emph{J. Comput. Appl. Math.}, 259: 336-349, 2014.
	
	\bibitem{Arriojas2007}M. Arriojas, Y.Z. Hu, S.E. Monhammed, and G. Pap, A delayed Black and Scholes formula. \emph{Stoch. Anal. Appl.}, 25(2): 471-492, 2007.
	
	\bibitem{BL2003}N. Bansal, Z. Liu, Capacity, delay and mobility in wireless ad-hoc networks. \emph{IEEE INFOCOM 2003. Twenty-second Annual Joint Conference of the IEEE Computer and Communications Societies (IEEE Cat. No. 03CH37428)}, 2: 1553-1563, 2003.
	
	\bibitem{BS2014}A. Bensoussan, C.C. Siu, S.C.P. Yam, and H.L. Yang, A class of non-zero-sum stochastic differential investment and reinsurance games. \emph{Automatica J. IFAC}, 50(8): 2025-2037, 2014.
	
	\bibitem{CW2010}L. Chen, Z. Wu, Maximum principle for the stochastic optimal control problem with delay and application. \emph{Automatica J. IFAC}, 46(6): 1074-1080, 2010.
	
	\bibitem{CWY2012}L. Chen, Z. Wu, and Z.Y. Yu, Delayed stochastic linear‐quadratic control problem and related applications. \emph{J. Appl. Math.}, 2012(1): 835319, 2012.
	
    \bibitem{CZ2026}P. Chen, F. Zhang, Maximum principle for partial information non-zero sum stochastic differential games with mixed delays. \emph{Automatica J. IFAC}, 183: 112570, 2026.

	\bibitem{Hamadene1999}S. Hamadène, Nonzero sum linear-quadratic stochastic differential games and backward-forward equations. \emph{Stochastic Anal. Appl.}, 17(1): 117-130, 1999.
	
	\bibitem{KSH2000}T. Kailath, A.H. Sayed, and B. Hassibi, \emph{Linear Estimation}, Prentice-Hall, Upper Saddle River, NJ, USA, 2000.
	
    \bibitem{MXW2022}T.F. Ma, J.J. Xu, and H.X. Wang, LQ control of It\^{o} stochastic system with asymmetric information. \emph{J. Math. Anal. Appl.}, 512: 126165, 2022.

	\bibitem{Ma2022}T.F. Ma, J.J. Xu, and H.S. Zhang, Explicit solution to forward and backward stochastic differential equations with state delay and its application to optimal control. \emph{Control Theory Technol.}, 20: 303-315, 2022.
	
	\bibitem{Ma2024}T.F. Ma, J.J. Xu, and H.S. Zhang, Explicit solution to delayed forward and backward stochastic differential equations. \emph{Internat. J. Systems Sci.}, 55(10): 2144-2153, 2024.
	
	\bibitem{Meng2025}W.J. Meng, J.T. Shi, T.X. Wang, and J.F. Zhang, A general maximum principle for optimal control of stochastic differential delay systems. \emph{SIAM J. Control Optim.}, 63(1): 175-205, 2025.
	
	\bibitem{NWY2022}T.Y. Nie, F.L. Wang, and Z.Y. Yu, Maximum principle for general partial information nonzero sum stochastic differential games and applications. \emph{Dyn. Games Appl.}, 12(2): 608-631, 2022.
	
	\bibitem{PY2009}S.G. Peng, Z. Yang, Anticipated backward stochastic differential equations. \emph{Ann. Probab.}, 37(3): 877-902, 2009.
	
	\bibitem{SW2015}J.T. Shi, G.C. Wang, A nonzero sum differential game of BSDE with time-delayed generator and applications. \emph{IEEE Trans. Automat. Control}, 61(7): 1959-1964, 2015.
	
	\bibitem{SWX2016}J.T. Shi, G.C. Wang, and J. Xiong, Leader-follower stochastic differential game with asymmetric information and applications. \emph{Automatica J. IFAC}, 63: 60-73, 2016.
	
	\bibitem{SWX2017}J.T. Shi, G.C. Wang, and J. Xiong, Linear-quadratic stochastic Stackelberg differential game with asymmetric information. \emph{Sci. China Inf. Sci.}, 60(9): 092202, 2017.
	
	\bibitem{Starr1969}A.W. Starr, Y.C. Ho, Nonzero-sum differential games. \emph{J. Optim. Theory Appl.}, 3(3): 184-206, 1969.

	\bibitem{Sun2019}J.R. Sun, J.M. Yong, Linear-quadratic stochastic two-person nonzero-sum differential games: open-loop and closed-loop Nash equilibria. \emph{Stoch. Process. Appl.}, 129(2): 381-418, 2019.
	
	\bibitem{WXX2018}G.C. Wang, H. Xiao, and J. Xiong, A kind of LQ non-zero sum differential game of backward stochastic differential equation with asymmetric information. \emph{Automatica J. IFAC}, 97: 346-352, 2018.
	
	\bibitem{WY2010}G.C. Wang, Z.Y. Yu, A Pontryagin's maximum principle for non-zero sum differential games of BSDEs with applications. \emph{IEEE Trans. Automat. Control}, 55(7): 1742-1747, 2010.
	
	\bibitem{WY2012}G.C. Wang, Z.Y. Yu, A partial information non-zero sum differential game of backward stochastic differential equations with applications, \emph{Automatica J. IFAC}, 48(2): 342-352, 2012.
	
	\bibitem{WZX2021}H.X. Wang, H.S. Zhang, and L.H. Xie, Optimal control and stabilization for It\^{o} systems with input delay. \emph{J. Syst. Sci. Complex.}, 34(5): 1895-1926, 2021.
	
	\bibitem{Wu2005}Z. Wu, Forward-backward stochastic differential equations, linear quadratic stochastic optimal control and nonzero sum differential games. \emph{J. Syst. Sci. Complex.}, 18(2): 179-192, 2005.
	
	\bibitem{WY2005}Z. Wu, Z.Y. Yu, Linear quadratic nonzero-sum differential games with random jumps. \emph{Appl. Math. Mech.}, 26(8): 1034-1039, 2005.
	
	\bibitem{XZ2018}J.J. Xu, H.S. Zhang, Solution to delayed FBSDEs and application to the stochastic LQ problem with input delay. \emph{IEEE Trans. Circuits Syst. II-Express Briefs}, 65(6): 769-773, 2018.
	
	\bibitem{XZX2018}J.J. Xu, H.S. Zhang, and L.H. Xie, General linear forward and backward stochastic difference equations with applications. \emph{Automatica J. IFAC}, 96: 40-50, 2018.
	
	\bibitem{YZ1999}J.M. Yong, X.Y. Zhou, \emph{Stochastic Controls: Hamiltonian Systems and HJB Equations}, Springer-Verlag, New York, 1999.
	
	\bibitem{ZX2017}H.S. Zhang, J.J. Xu, Control for It\^{o} stochastic systems with input delay. \emph{IEEE Trans. Automat. Control}, 62(1): 350-365, 2017. 

\end{thebibliography}
\end{document}